# Scattering theory for the Hodge Laplacian

Robert Baumgarth[1]

[1]Department of Mathematics, Université du Luxembourg, Maison du Nombre, 6, Avenue de la Fonte, 4364 Esch-sur-Alzette, Grand-Duché de Luxembourg

robert.baumgarth@uni.lu

4th June 2021

**Abstract**

We prove using an integral criterion the existence and completeness of the wave operators $W_{\pm}(\Delta_h^{(k)}, \Delta_g^{(k)}, I_{g,h}^{(k)})$ corresponding to the Hodge Laplacians $\Delta_\nu^{(k)}$ acting on differential $k$-forms, for $\nu \in \{g, h\}$, induced by two quasi-isometric Riemannian metrics $g$ and $h$ on a complete open smooth manifold $M$. In particular, this result provides a criterion for the absolutely continuous spectra $\sigma_{\mathrm{ac}}(\Delta_g^{(k)}) = \sigma_{\mathrm{ac}}(\Delta_h^{(k)})$ of $\Delta_\nu^{(k)}$ to coincide. The proof is based on gradient estimates obtained by probabilistic Bismut-type formulae for the heat semigroup defined by spectral calculus. By these localised formulae, the integral criterion requires local curvature bounds and some upper local control on the heat kernel acting on functions provided the Weitzenböck curvature endomorphism is in the Kato class, but no control on the injectivity radii. A consequence is a stability result of the absolutely continuous spectrum under a Ricci flow. As an application we concentrate on the important case of conformal perturbations.

**Contents** 1

## 1 Introduction

Let $(M, g)$ be a non-compact geodesically complete Riemannian manifold without boundary. The Hodge Laplacian $\Delta_g^{(k)}$ acting on differential $k$-forms carries important geometric and topological information about $M$, of particular interest is the spectrum $\sigma(\Delta_g^{(k)})$ of $\Delta_g^{(k)}$. If $M$ is compact, then the spectrum consists of eigenvalues with finite multiplicity. If $M$ is non-compact, then the spectrum contains some absolutely continuous part (cf. [RS79; Wei80]). A natural question to ask is to what extent can we control the absolutely continuous part of $\sigma(\Delta_g^{(k)})$ and under which assumptions on the geometry of $(M, g)$?

A systematic approach to control the absolutely continuous part of the spectrum $\sigma_{\mathrm{ac}}(\Delta_g^{(k)})$ is inspired by quantum mechanics, namely scattering theory: Assume that there is another Riemannian metric $h$ on $M$ such that $h$ is quasi-isometric to $g$, i.e. there exists a constant $C \geqslant 1$ such that $(1/C)g \leqslant h \leqslant Cg$. We show that under suitable assumptions the wave operators

$$W_\pm(\Delta_h^{(k)}, \Delta_g^{(k)}, I_{g,h}^{(k)}) = \operatorname*{s-lim}_{t\to\pm\infty} \mathrm{e}^{it\Delta_h^{(k)}} I_{g,h}^{(k)} \mathrm{e}^{-it\Delta_g^{(k)}} \mathrm{P}_{\mathrm{ac}}(\Delta_g^{(k)})$$

exist and are complete, where the limit is taken in the strong sense, and

$$I_{g,h}^{(k)} : \Omega^k_{\mathsf{L}^2}(M, g) \to \Omega^k_{\mathsf{L}^2}(M, h)$$

denotes a bounded identification operator between the Hilbert spaces of equivalence classes of square-integrable Borel $k$-forms on $M$ corresponding to the metric $g$ and $h$ respectively (cf. Theorem A.3 in the Appendix and Section 2 for details). Then as well-known (cf. appendix A, it follows in particular that

$$\sigma_{\mathrm{ac}}(\Delta_h^{(k)}) = \sigma_{\mathrm{ac}}(\Delta_g^{(k)}).$$

Considering Laplacians acting on 0-forms, i.e. functions, on $M$, Müller & Salomonsen [MS07] studied the existence and completeness of the wave operators corresponding to the Laplace-Beltrami operator by assuming both metrics to have a $C^\infty$-bounded geometry and a weighted integral condition involving a second order deviation of the metrics. Hempel, Post & Weder [HPW14] improved the result of [MS07] by assuming only a zeroth order deviation of the metric $g$ from $h$ and a weighted integral condition involving a local lower bound of the injectivity radius and the Ricci curvatures. However, detailed control on the sectional curvature is needed to get control over the injectivity radii. In general, injectivity radii are hard to calculate.

Recently, Güneysu & Thalmaier [GT20] established a rather simple integral criterion induced by two quasi-isometric Riemannian metrics only depending on a local upper bound on the heat kernel and certain explicitly given local lower bound on the Ricci curvature using stochastic methods, namely a Bismut-type formula for the derivative of the heat semigroup [AT10].

Considering Laplacians acting on differential $k$-forms, Bei, Güneysu & Müller [BGM17] generalised the previous results in [MS07] for the case of *conformally equivalent metrics* under a mild first order control on the conformal factor.

Using a similar method, very recently Boldt & Güneysu [BG20] extended the result of [GT20] to a non-compact spin manifold with a fixed topological spin structure and two complete Riemannian metrics with bounded sectional curvatures. As the metrics induce Dirac operators $\mathsf{D}_g$ and $\mathsf{D}_h$, they can show existence and completeness of the wave operators corresponding to the Dirac operators $W_\pm(\mathsf{D}_h, \mathsf{D}_g, I_{g,h})$ and their squares $W_\pm(\mathsf{D}_h^2, \mathsf{D}_g^2, I_{g,h})$.

*In this paper, we address the natural question: Can we extend the result of [GT20] to the setting of differential $k$-forms for two quasi-isometric Riemannian metrics?*

We will show that this can be done if the Weitzenböck curvature endomorphisms is in the Kato class and assuming an integral criterion only depending on a local upper bound on the heat kernel and certain explicitly given local curvature bounds. In addition, a necessary assumption will be a bound on a weight function measuring the first order deviation of the metrics in terms of the corresponding covariant derivatives $\nabla^h$ and $\nabla^g$ which is a one-form on $M$ with values in $\operatorname{End} \mathrm{T}M$.

Therefore, we consider the Hodge Laplacian $\Delta_g$, also known as Laplace-de Rham operator, acting on the full exterior bundle $\Omega(M) = \Gamma(\bigwedge \mathrm{T}^*M)$, i.e. the complex separable Hilbert space of differential forms on $M$. The Hodge Laplacian $\Delta_g$ is related to the horizontal Laplacian $\square_g = -(\nabla^g)^*\nabla^g$ by the Weitzenböck formula $\Delta_g = \square_g - \mathscr{R}_g$, where Weitzenböck curvature operator $\mathscr{R}_g \in \Gamma(\operatorname{End} \Omega(M))$ is a symmetric field of endomorphisms. In particular, when acting on 1-forms, $\mathscr{R}_g^{\mathrm{tr}}|_{\Omega^1(M,g)} = \mathrm{Ric}_g$ and, acting on functions, $\mathscr{R}_g^{\mathrm{tr}}|_{\Omega^0(M,g)} = 0$. We assume that $\mathscr{R}_g$ is in the Kato class, i.e. that the fibrewise taken operator norm $|\mathscr{R}_g|_g$ (which is a Borel function on $M$) of $\mathscr{R}_g$ is in the Kato class (cf. Definition 3.3). We are now in the position to state our main result, cf. Theorem 4.3 below.

**Main result.** *Assume that $g$ and $h$ are two geodesically complete and quasi-isometric Riemannian metrics on $M$, denoted $g \sim h$, and assume that there exists $C < \infty$ such that $\left|\delta^{\nabla}_{g,h}\right| \leqslant C$, and that for both $\nu \in \{g, h\}$, $\mathscr{R}_\nu$ is in the Kato class and it holds*

$$\int \max\left\{\delta_{g,h}(x), \delta^{\nabla}_{g,h}(x) + \Xi_g(x,s), \Psi_\nu(x,s)\right\} \Phi_\nu(x,s)\,\mathrm{vol}_\nu(\mathrm{d}x) < \infty, \qquad \textit{some } s > 0, \tag{1.1}$$

*where*

- $\mathrm{vol}_\nu$ *denotes the Riemannian volume measure with respect to the metric $\nu$,*
- $\Psi_\nu(x,s) : M \to (0,\infty)$ *is a function explicitly given terms of local curvature bounds (cf. (3.14) in Section 3) and a finite constant $c_\gamma(\underline{\mathscr{R}}^-)$ (cf. (3.5) in Section 3),*
- $\Xi_\nu(\cdot,s) : M \to (0,\infty)$ *is a function explicitly given in terms of $\Psi_\nu(x,s)$ and a local bound on the derivative of the curvature (cf. (3.20)),*
- $\Phi_\nu(\cdot,s) : M \to (0,\infty)$ *is a local upper bound on the heat kernel acting on functions on $(M,\nu)$ (cf. (2.15) in Section 2),*
- $\delta_{g,h} : M \to (0,\infty)$ *a zeroth order deviation of the metrics from each other (cf. (2.9) in Sec. 2),*
- $\delta^{\nabla}_{g,h} : M \to [0,\infty)$ *a first order deviation of the metrics (cf. (2.10) in Section 2).*

*Then the wave operators $W_\pm(\Delta_h, \Delta_g, I_{g,h})$ exist and are complete. Moreover, $W_\pm(\Delta_h, \Delta_g, I_{g,h})$ are partial isometries with initial space $\operatorname{ran} \mathrm{P}_{\mathrm{ac}}(\Delta_g)$ and final space $\operatorname{ran} \mathrm{P}_{\mathrm{ac}}(\Delta_h)$. In particular,*

$$\sigma_{\mathrm{ac}}(\Delta_g) = \sigma_{\mathrm{ac}}(\Delta_h).$$

We will see that a zeroth order deviation $\delta_{g,h}$ of the metrics from each other is induced by quasi-isometry. In comparison to the case of 0-forms, i.e. functions, it turns out, working on higher degree differential forms, also a first order deviation of the metrics $\delta^{\nabla}_{g,h} = \left|\nabla^h - \nabla^g\right|_g^2$ is necessary. But note that $\nabla^h - \nabla^g$ is a one-form on $M$ with values in $\operatorname{End} \mathrm{T}M$.

In contrast to previous results, it seems that we are the first to assume global curvature conditions in terms of the Kato class, more precisely, that the Weitzenböck curvature endomorphism is in the Kato class.

To this end, our strategy is to verify the assumptions given by a variant of the Belopol'skii-Birman theorem A.3 which is adapted to our special case of two Hilbert spaces. The main technical difficulty is to show that the operator

$$T_s^{g,h} = \Delta_h^{(k)} e^{s\Delta_h^{(k)}} I_{g,h} e^{s\Delta_g^{(k)}} - e^{s\Delta_h^{(k)}} I_{g,h} e^{s\Delta_g^{(k)}} \Delta_g^{(k)}$$

is trace class. As the product of Hilbert-Schmidt operators is trace class, our idea is to decompose the operator $T_s^{g,h}$ in such a way that the terms only consist of (transformed) derivations of Hilbert-Schmidt estimates and bounded multiplication operators. In comparison to thet corresponding decomposition formula in [GT20, Lemma 4.1], the analysis becomes considerably more difficult because the quadratic form associated to $\Delta_h^{(k)}$ involves not only the exterior derivative $\mathbf{d}^{(k)}$, but also the codifferential $\boldsymbol{\delta}_h^{(k)}$ which depends on the metric by definition. Moreover, we encounter quantities transformed by a smooth vector bundle morphism $\mathscr{A}_{g,h}$ induced by the quasi-isometry (cf. (2.4) below). Using the quasi-isometry of the metrics we can give a formula how to express the codifferential $\boldsymbol{\delta}_h$ with respect the metric $h$ in terms of the codifferential $\boldsymbol{\delta}_g$ in terms of $g$ (cf. Lemma 2.8). Using the metric description for the exterior derivative and the codifferential (cf. Lemma 3.11), we can express the corresponding quantities transformed by $\mathscr{A}_{g,h}$ in terms of the covariant derivative $\nabla^g$ of $g$ applied to the semigroup (cf. Proposition 3.15).

Our tool to obtain the Hilbert-Schmidt estimates for various derivatives of the heat semigroup will be derived probabilistic Bismut-type derivative formulae for the exterior derivative, codifferential and covariant derivative (cf. Theorem 3.9 and 3.10) following the ideas in [DT01]. The gradient estimates are then a direct consequence and the probabilistic formulae used provide us, in particular, with explicit local constants.

We first look at total differential forms, then finally everything filters through the form degree to differential $k$-forms. The particularly important case, in which two quasi-isometric Riemannian metrics differ by a conformal metric change, is a direct consequence of our main result.

Because our result is independent of the injectivity radii we have the following application to the Ricci flow. Let $\mathrm{R}_g$ be the Riemannian curvature tensor with respect to the metric $g$ and set $\dim M =: m$.

**Corollary 6.1.** *Let $S > 0$, $\lambda \in \mathbb{R}$ and assume that*

*(a) the family $(g_s)_{0\leqslant s\leqslant S} \subset \mathrm{Metr}M$ evolves under a Ricci-type flow*

$$\partial_s g_s = \lambda \,\mathrm{Ric}_{g_s}, \qquad \text{for all } \ 0 \leqslant s \leqslant S$$

*(b) the initial metric $g_0$ is geodesically complete*

*(c) there is some $C > 0$ such that* $\left|\mathrm{R}_{g_s}\right|_{g_s}, \left|\nabla^{g_s}\mathrm{R}_{g_s}\right|_{g_s} \leqslant C \qquad \text{for all } \ 0 \leqslant s \leqslant S.$

*We set, for all $x \in M$,*

$$M_1(x) := \sup \left\{ \left| \mathrm{Ric}_{g_s}(v,v) \right| : 0 \leqslant s \leqslant S,\ v \in \mathsf{T}_x M,\ |v|_{g_s} \leqslant 1 \right\},$$
$$M_2(x) := \sup \left\{ \left| \nabla^{g_s}_v \mathrm{Ric}_{g_s}(u,w) + \nabla^{g_s}_u \mathrm{Ric}_{g_s}(v,w) + \nabla^{g_s}_w \mathrm{Ric}_{g_s}(u,v) \right| : 0 \leqslant s \leqslant S,\right.$$
$$\left. u,v,w \in \mathsf{T}_x M,\ |u|_{g_s}, |v|_{g_s}, |w|_{g_s} \leqslant 1 \right\}.$$

*Let $\mathsf{B}_g(x,R)$ denote the open geodesic ball. If*

$$\int \mathrm{vol}_{g_0}(\mathsf{B}_{g_0}(x,1))^{-1} \max \left\{ \sinh\left( \frac{m}{4} S\, |\lambda|\, M_1(x) \right), M_2(x) \right\} \mathrm{vol}_{g_0}(\mathrm{d}x) < \infty,$$

*then $\sigma_{\mathrm{ac}}(\Delta_{g_s}) = \sigma_{\mathrm{ac}}(\Delta_{g_0})$ for all $0 \leqslant s \leqslant S$.*

Thereupon, we reify our main results to the case of global curvature bounds: The curvature operator (with respect to the metric $g$) $Q_g$ is uniquely determined by the equation

$$\big(Q_g(X \wedge Y), U \wedge V\big)_g = \big(\mathsf{R}_g(X,Y)U, V\big)_g$$

for all smooth vector fields $X,Y,U,V \in \Gamma_{C^\infty}(\mathsf{T}M)$. By the Gallot-Meyer estimate [GM75], a global bound $Q_g \geqslant -K$, for some constant $K > 0$, already implies that the curvature endomorphism in the Weitzenböck formula (3.1) is globally bounded by $\mathscr{R}_g^{(k)} \geqslant -Kk(m-k)$.

Then the function $\Xi_g(x,s)$ simplifies to $\Theta_g(x,s) := \left(1 + \max_{y \in \mathsf{B}_g(x,1)} \left|\nabla^g \mathsf{R}_g(y)\right|\right)^2$. In this case, our main result reads as follows.

**Theorem 6.8.** *Let $Q_\nu \geqslant -K$, for some constant $K > 0$ for both $\nu \in \{g,h\}$. Let $g,h \in \mathrm{Metr}M$ such that $g \sim h$ and assume that there exists $C < \infty$ such that $\left|\delta^\nabla_{g,h}\right| \leqslant C$ and that for some (then both by quasi-isometry) $\nu \in \{g,h\}$*

$$\int \max \left\{ \delta_{g,h}(x), \delta^\nabla_{g,h}(x) + \Theta_g(x,s) \right\} \Phi_\nu(x,s)\, \mathrm{vol}_\nu(\mathrm{d}x) < \infty, \qquad \textit{some } s > 0.$$

*Then, the wave operators $W_\pm(\Delta_h, \Delta_g, I)$ exist and are complete. Moreover, $W_\pm(\Delta_h, \Delta_g, I)$ are partial isometries with initial space $\operatorname{ran} \mathsf{P}_{\mathrm{ac}}(\Delta_g)$ and final space $\operatorname{ran} \mathsf{P}_{\mathrm{ac}}(\Delta_h)$, and we have $\sigma_{\mathrm{ac}}(\Delta_g) = \sigma_{\mathrm{ac}}(\Delta_h)$.*

A direct consequence and additional application of our main result is the particularly important case of conformal perturbations under local and global curvature bounds, generalising the results in [BGM17].

A final application is provided through a result by Cheeger, Fukaya and Gromov [CFG92] known as Cheeger-Gromov's thick/thin decomposition: On any complete Riemannian $m$-manifold $(M,g)$ with bounded *sectional curvature* $\left|\kappa_g\right| \leqslant 1$, there exists a Riemannian metric $g_\varepsilon$ on $M$ such that $g_\varepsilon$ is $\varepsilon$-quasi-isometric to $g$ and has bounded covariant derivatives. Hence, in this case, the assumptions of our main result may be suitably relaxed, cf. Theorem 6.11.

Let us end the introduction with a short outline of the paper. Section 2 introduces the necessary notation and deviation maps. In Section 3, we calculate bounds for exterior derivative, codifferential and covariant derivative of the heat semigroup defined by spectral calculus using

Bismut-type derivative formulae. Our main results are explained in Section 4. After this, we prove the main result in Section 5 by making use of a variant of the abstract Belopol'skii-Birman Theorem A.3.1. We close in Section 6 with applications to the Ricci flow 6.1, state the main result in case of differential $k$-forms 6.2, the particularly important cases of conformal perturbations 6.3, specify our results for global curvature bounds 6.4 and $\varepsilon$-close Riemannian metrics 6.5.

## 2 Setting & Notation

Let $(M, g)$ be a complete smooth Riemannian manifold without boundary of dimension $m := \dim M \geqslant 2$ and $(\cdot,\cdot)_g$ its Riemannian metric. We write $\mathrm{vol}_g$ for the corresponding volume measure (with respect to the metric $g$) and denote by $\mathsf{Metr}M$ the set of all smooth Riemannian metrics on $M$. All bundles will be understood complexified, e.g. the full exterior bundle

$$\bigwedge \mathsf{T}^*M = \bigoplus_{k=0}^{m} \bigwedge^k \mathsf{T}^*M, \quad \text{with the usual convention } \bigwedge^0 \mathsf{T}^*M = \mathbb{C}. \tag{2.1}$$

Given smooth complex vector bundles $E_1 \to M$ and $E_2 \to M$ the complex linear space of smooth linear partial differential operators from $E_1$ to $E_2$ of order $\leqslant k \in \mathbb{N}_0$ is denoted by $\mathscr{D}^{(k)}(M; E_1, E_2)$, with shorthand notation $\mathscr{D}^{(k)}(M; E_1)$ if $E_1 \equiv E_2$. On a vector bundle $E \to M$ (e.g. $E = \bigwedge^k \mathsf{T}^*M$) the corresponding fibre norms are denoted by

$$|\varphi|_g := (\varphi, \varphi)_g^{1/2} \qquad \text{for any section } \varphi \in \Gamma(E),$$

where $\Gamma(E) := \Gamma_{\mathsf{C}^\infty}(E)$ denotes all *smooth* sections of $E$ and $\Gamma_{\mathsf{L}^2}(E)$ the $\mathsf{L}^2$-section of $E$.

In the case of $E = \bigwedge^k \mathsf{T}^*M$, we indicate the corresponding form degree by an index: For example, $\nabla^{g,(k)}$ or $(\cdot,\cdot)_g^{(k)}$ etc.

We denote by $\Omega_{\mathsf{L}^2}(M, g) := \Gamma_{\mathsf{L}^2}(\bigwedge \mathsf{T}^*M)$ the complex separable Hilbert space of equivalence classes $\alpha$ of square-integrable Borel forms on $M$ such that

$$\|\alpha\|_g^2 := \|\alpha\|^2_{\Omega_{\mathsf{L}^2}(M,g)} := \int_M |\alpha(x)|_g^2 \,\mathrm{vol}_g(\mathrm{d}x) < \infty,$$

with inner product

$$\langle \alpha, \beta \rangle_g := \langle \alpha, \beta \rangle_{\Omega_{\mathsf{L}^2}(M,g)} := \int_M (\alpha(x), \beta(x))_g \,\mathrm{vol}_g(\mathrm{d}x).$$

Analogously, we write $\Omega^k_{\mathsf{L}^2}(M, g)$ for the Hilbert space of Borel $k$-forms. In particular,

$$\Omega_{\mathsf{L}^2}(M, g) = \bigoplus_{k=0}^{m} \Omega^k_{\mathsf{L}^2}(M, g).$$

To relax notation, we set

$$\Omega(M, g) := \Omega_{\mathsf{C}^\infty}(M, g) \qquad \text{and} \qquad \Omega^k(M, g) := \Omega^k_{\mathsf{C}^\infty}(M, g).$$

for the set of all *smooth* forms, and *smooth* $k$-forms respectively, on $(M, g)$.

Further, for some $\alpha \in \Omega^k(M)$, we denote by

$$\bullet \wedge \alpha \in \mathscr{D}^{(0)}(M; \bigwedge^k \mathsf{T}^*M, \bigwedge^{k+1} \mathsf{T}^*M)$$

the exterior product and its formal adjoint with respect to $g$, the interior multiplication, by

$$\bullet \lrcorner_g \alpha := (\bullet \wedge \alpha)^{*_g} \in \mathscr{D}^{(0)}(M; \textstyle\bigwedge^k \mathsf{T}^*M, \bigwedge^{k-1} \mathsf{T}^*M).$$

The interior multiplication corresponds to the contraction of $\alpha \in \Omega^k(M)$ with a vector field $X \in \boldsymbol{\Gamma}(\mathsf{T}M)$, i.e.

$$X \lrcorner_g \alpha(X_1, \dots, X_{k-1}) := \alpha(X, X_1, \dots, X_{k-1}), \qquad \text{for all } X_1, \dots, X_{k-1} \in \boldsymbol{\Gamma}(\mathsf{T}M),$$

and is an antiderivation, to wit

$$X \lrcorner_g (\alpha \wedge \beta) = (X \lrcorner \alpha) \wedge \beta + (-1)^k \alpha \wedge (X \lrcorner \beta) \qquad \text{for all } \alpha \in \Omega^k(M),\ \beta \in \Omega(M).$$

We denote by

$$\mathbf{d}^{(k)} \in \mathscr{D}^{(1)}(M; \textstyle\bigwedge^k \mathsf{T}^*M, \bigwedge^{k+1} \mathsf{T}^*M)$$
$$\boldsymbol{\delta}_g^{(k)} \in \mathscr{D}^{(1)}(M; \textstyle\bigwedge^k \mathsf{T}^*M, \bigwedge^{k-1} \mathsf{T}^*M)$$

the exterior derivative on $k$-forms and, respectively, the codifferential as the formal adjoint of $\mathbf{d}^{(k-1)}$. Then the Hodge Laplacian can be written as the sum

$$\Delta_g^{(k)} := -\left(\boldsymbol{\delta}_g^{(k+1)}\mathbf{d}^{(k)} + \mathbf{d}^{(k-1)}\boldsymbol{\delta}_g^{(k)}\right) \in \mathscr{D}^{(2)}(M; \textstyle\bigwedge^k \mathsf{T}^*M)$$

and its Friedrichs realisation in $\Omega^k_{\mathsf{L}^2}(M, g)$ will be again denoted by $-\Delta_g^{(k)} \geqslant 0$. In particular, for $k = 0$, we recover the special case of the Laplace-Beltrami operator acting on 0-forms, i.e. functions,

$$\Delta_g^{(0)} = -\boldsymbol{\delta}_g^{(1)}\mathbf{d}^{(0)} \in \mathscr{D}^{(2)}(M).$$

Furthermore, we set

$$\mathbf{d} := \bigoplus_{k=0}^{m} \mathbf{d}^{(k)} \in \mathscr{D}^{(1)}(M; \textstyle\bigwedge \mathsf{T}^*M)$$
$$\boldsymbol{\delta}_g := \bigoplus_{k=0}^{m} \boldsymbol{\delta}_g^{(k)} \in \mathscr{D}^{(1)}(M; \textstyle\bigwedge \mathsf{T}^*M)$$

and define the underlying Dirac-type operator $\mathbf{D}_g$, and the (total) Hodge Laplacian $\Delta_g$

$$\mathbf{D}_g := \mathbf{d} + \boldsymbol{\delta}_g \in \mathscr{D}^{(1)}(M; \textstyle\bigwedge \mathsf{T}^*M)$$
$$\Delta_g := -\mathbf{D}_g^2 \in \mathscr{D}^{(2)}(M; \textstyle\bigwedge \mathsf{T}^*M)$$

where the Friedrichs realisation of $\Delta_g$ in $\Omega^k_{\mathsf{L}^2}(M, g)$ will be again denoted by $-\Delta_g \geqslant 0$. In particular,

$$\Delta_g\big|_{\Omega^k(M,g)} = \Delta_g^{(k)} \in \mathscr{D}^{(2)}(M; \textstyle\bigwedge^k \mathsf{T}^*M)$$

and

$$\Delta_g = \bigoplus_{k=0}^{m} \Delta_g^{(k)} \qquad \text{as self-adjoint operators.}$$

Since $g$ is (geodesically) complete, it follows that the operators $\mathbf{D}_g, \Delta_g, \Delta_g^{(k)}$ are essentially self-adjoint on the corresponding space of smooth compactly supported forms [Str83].

Next, recall that for the $k$-fold exterior product of the vector space $\mathrm{T}^*M$, we obtain a scalar product $(\cdot,\cdot)_g$ on $\bigwedge^k \mathrm{T}^*M$ by the bilinear extension of

$$\left(\alpha_1 \wedge \ldots \wedge \alpha_k, \beta_1 \wedge \ldots \wedge \beta_k\right)_g = \det\left(\alpha_k, \beta_l\right)_g . \tag{2.2}$$

Any $A \in \mathrm{End}(\mathrm{T}^*M)$ induces a linear map

$$\begin{aligned}\textstyle\bigwedge^m A : \bigwedge^m \mathrm{T}^*M &\to \textstyle\bigwedge^m \mathrm{T}^*M\\ \alpha_1 \wedge \ldots \wedge \alpha_m &\mapsto A\alpha_1 \wedge \ldots \wedge A\alpha_m.\end{aligned}$$

As $\bigwedge^m \mathrm{T}^*M$ is one-dimensional, the map $\bigwedge^m A$ is given by multiplication with a unique number, denoted by $\det A$,

$$\textstyle\bigwedge^m A(e_1 \wedge \ldots \wedge e_m) = (\det A)\, e_1 \wedge \ldots \wedge e_m,$$

where $(e_1, \ldots, e_m)$ is a basis for $\bigwedge^m \mathrm{T}_x^*M$.

A Riemannian metric $(u,v)_g = g(u,v)$ for $u, v \in \mathrm{T}_xM$ gives by definition an inner product on each tangent space $\mathrm{T}_xM$ ($x \in M$). By Riesz' representation theorem, $g$ provides a natural isomorphism between tangent and cotangent bundle given by $v \mapsto (v,\cdot)_g$,

$$\mathrm{T}M \underset{\sharp^g}{\overset{\flat^g}{\rightleftarrows}} \mathrm{T}^*M.$$

More precisely, we define the **sharp operator** $\sharp^g$ (with respect to $g$) by

$$\sharp^g : \mathrm{T}^*M \to \mathrm{T}M, \qquad \alpha(v) = g(\alpha^{\sharp^g}, v).$$

The Riemannian metric $g$ defines a metric $g$ on $\mathrm{T}^*M$, the **cometric**, via

$$g(\alpha,\beta) := g(\alpha^{\sharp^g}, \beta^{\sharp^g}) \qquad \text{for all } \ \alpha, \beta \in \mathrm{T}_x^*M,\ x \in M,$$

which extends to a metric on $\bigwedge^k \mathrm{T}^*M$ according to (2.2).

**Definition 2.1.** A smooth Riemannian metric $h \in \mathrm{Metr}M$ is called **quasi-isometric to** $g$, denoted $g \sim h$, if there exists a constant $C \geqslant 1$ such that (to be understood pointwise, as bilinear forms)

$$\frac{1}{C} g \leqslant h \leqslant C g.$$

Given $g, h \in \mathrm{Metr}M$, we define a vector bundle morphism

$$A := A_{g,h} : \mathrm{T}M \xrightarrow{\sim} \mathrm{T}M, \qquad h(u,v) = g(Au, v), \qquad \text{for all } \ x \in M,\ u, v \in \mathrm{T}_xM. \tag{2.3}$$

Note that the vector bundle morphism $A = A_{g,h}$ induces a vector bundle morphism on the cotangent bundle via

$$A : \mathrm{T}^*M \xrightarrow{\sim} \mathrm{T}^*M, \qquad \alpha \mapsto A\alpha := \alpha \circ A.$$

**Lemma 2.2** (and **Definition**). *In terms of the notations above, we have*

$$h(\alpha,\beta) = g(A^{-1}\alpha,\beta) \qquad \text{for all } \alpha,\beta \in \mathsf{T}_x^*M,\ x\in M.$$

*Extending* $A^{-1} = A_{g,h}^{-1}$ *to a smooth vector bundle morphism by*

$$\mathscr{A} := \mathscr{A}_{g,h}(x) := (\bigwedge A_{g,h}^{-1})_x : \bigwedge \mathsf{T}^*M \xrightarrow{\sim} \bigwedge \mathsf{T}^*M, \qquad \mathscr{A}\alpha := \alpha\circ\mathscr{A} \tag{2.4}$$

*we obtain*

$$g(\mathscr{A}_{g,h}(x)\alpha,\beta) = h(\alpha,\beta) \qquad \text{for } x\in M,\ \alpha,\beta\in\bigwedge \mathsf{T}_x^*M.$$

In the following the induced metrics will be understood complexified (conjugate-linear in the first variable and linear in the second).

**Remark 2.3.** (i) By the positive-definiteness of $h$ (or $g$), $\mathscr{A}_{g,h}(x)$ has only positive eigenvalues ($x\in M$). By the symmetry of $g$ and $h$ the endomorphism $\mathscr{A}_{g,h}$ is fibrewise self-adjoint with respect to $g$ and $h$. Therefore, the fibrewise operator norm $|\cdot|_g$ (or $|\cdot|_h$) induced by the metric $g$ (or $h$) of $\mathscr{A}$ is equivalent to absolute value of the largest eigenvalue on the given fibre for both $g$ and $h$. Thus to relax notation, we may suppress the metric and simply write $|\mathscr{A}|$.

(ii) By the very definition, $\mathscr{A}^{1/2}$ is a (pointwise) isometry from $(\bigwedge \mathsf{T}^*M, g)$ to $(\bigwedge \mathsf{T}^*M, h)$.

**Proof of Lemma 2.2.** We prove the Lemma in several steps.

1° We calculate the sharp operator in the new metric. For $x\in M$, let $v\in\mathsf{T}_xM$ and $\alpha\in\mathsf{T}_x^*M$. By duality,

$$g\left(\alpha^{\sharp g}, v\right) = \alpha(v) = h\left(\alpha^{\sharp h}, v\right) = g\left(A\alpha^{\sharp h}, v\right) \qquad \Longrightarrow \qquad \alpha^{\sharp g} = A\alpha^{\sharp h}, \tag{2.5}$$

for all $v\in\mathsf{T}_xM$.

2° Let $\alpha,\beta\in\mathsf{T}_x^*M$, then

$$\begin{aligned} h(\alpha,\beta) &= h\left(\alpha^{\sharp^h}, \beta^{\sharp^h}\right) \\ &= g\left(A\alpha^{\sharp^h}, \beta^{\sharp^h}\right) = g(\alpha^{\sharp^h}, A\beta^{\sharp^h}) = g(A^{-1}\alpha^{\sharp^g}, \beta^{\sharp^g}) = g(\alpha\circ A^{-1},\beta), \end{aligned}$$

where we used that $A^{-1}\alpha^{\sharp g} = (\alpha\circ A^{-1})^{\sharp g}$ in the last equality.

3° For any $\alpha,\beta\in\Omega^k(M)$,

$$\begin{aligned} h(\alpha,\beta) &= h\left(\alpha_1^{\sharp h}\wedge\ldots\wedge\alpha_k^{\sharp h}, \beta_1^{\sharp h}\wedge\ldots\wedge\beta_k^{\sharp h}\right) \\ &\overset{(2.2)}{=} \det\left(\alpha_k^{\sharp h}, \beta_l^{\sharp h}\right)_h \\ &= \det\left(A\alpha_k^{\sharp h}, \beta_l^{\sharp h}\right)_g \\ &\overset{(2.5)}{=} \det\left(\alpha_k^{\sharp g}, A^{-1}\beta_l^{\sharp g}\right)_g \\ &\underset{\text{s.a.}}{\overset{A}{=}} \det\left(A^{-1}\alpha_k^{\sharp g}, \beta_l^{\sharp g}\right)_g \\ &= g\left(A^{-1}\alpha_1^{\sharp g}\wedge\ldots\wedge A^{-1}\alpha_k^{\sharp g}, \beta_1^{\sharp g}\wedge\ldots\wedge\beta_k^{\sharp g}\right) \\ &= g\left(\bigwedge^k A^{-1}\left(\alpha_1\wedge\ldots\wedge\alpha_k\right), \beta_1\wedge\ldots\wedge\beta_k\right) = g\left(\mathscr{A}\alpha,\beta\right) \end{aligned}$$

■

The following estimates are necessary tools for the main proof noting that they are independent of the quasi-isometry of $g$ and $h$.

**Lemma 2.4.** *Let $\mathscr{A} := \mathscr{A}_{g,h}$ be the smooth vector bundle morphism defined by (2.4). For any vector field $X \in \Gamma(TM)$, we get*

$$\nabla^g_X \mathscr{A} = \mathscr{A}(\nabla^h_X - \nabla^g_X) + (\nabla^h_X - \nabla^g_X)^* \mathscr{A},$$

*and the pointwise estimate*

$$\left|\nabla^g_X \mathscr{A}\right|_g \leqslant 2\,|\mathscr{A}|\,\left|\nabla^h_X - \nabla^g_X\right|_g, \tag{2.6}$$

*where $|\cdot|_g$ denotes the operator norm induced by the inner product $g$.*

**Proof.** We divide the proof into two steps.

**1°** If we differentiate the identity

$$(\alpha, \beta)_h = (\mathscr{A}\alpha, \beta)_g$$

in direction of $X$, we get on the one hand

$$\begin{aligned} X\,(\alpha, \beta)_h &= \left(\nabla^h_X \alpha, \beta\right)_h + \left(\alpha, \nabla^h_X \beta\right)_h \\ &= \left(\mathscr{A}\nabla^h_X \alpha, \beta\right)_g + \left(\mathscr{A}\alpha, \nabla^h_X \beta\right)_g \end{aligned}$$

and on the other hand

$$\begin{aligned} X\,(\mathscr{A}\alpha, \beta)_g &= \left(\nabla^g_X(\mathscr{A}\alpha), \beta\right)_g + \left(\mathscr{A}\alpha, \nabla^g_X \beta\right)_g \\ &= \left((\nabla^g_X \mathscr{A})\alpha, \beta\right)_g + \left(\mathscr{A}\nabla^g_X \alpha, \beta\right)_g + \left(\mathscr{A}\alpha, \nabla^g_X \beta\right)_g . \end{aligned}$$

Hence,

$$\left((\nabla^g_X \mathscr{A})\alpha, \beta\right)_g = \left(\mathscr{A}(\nabla^h_X - \nabla^g_X)\alpha, \beta\right)_g + \left(\alpha, \mathscr{A}(\nabla^h_X - \nabla^g_X)\beta\right)_g ,$$

using the self-adjointness of $\mathscr{A}$.

**2°** We estimate

$$\begin{aligned} \left|\nabla^g_X \mathscr{A}\right|_g &\leqslant \sup_{|v|,|w|\leqslant 1} \left|\left((\nabla^g_X \mathscr{A})v, w\right)_g\right| \\ &\leqslant \sup_{|v|,|w|\leqslant 1} \left(\left|\left(\mathscr{A}(\nabla^h_X - \nabla^g_X)v, w\right)_g\right| + \left|\left(v, \mathscr{A}(\nabla^h_X - \nabla^g_X)w\right)_g\right|\right) \\ &\leqslant 2\,|\mathscr{A}|\,\left|\nabla^h_X - \nabla^g_X\right|_g . \end{aligned}$$

∎

We point out that a similar argument was recently developed in [BG20] to prove the estimates in Lemma 2.5.

**Lemma 2.5.** *Let $\mathscr{A} := \mathscr{A}_{g,h}$ be the smooth vector bundle morphism defined by (2.4). For any vector field $X \in \Gamma(TM)$, we get the pointwise estimates*

$$\begin{aligned} \left|\nabla^g_X \mathscr{A}^{1/2}\right|_g &\leqslant |\mathscr{A}|\,\left|\mathscr{A}^{-1/2}\right|\,\left|\nabla^h_X - \nabla^g_X\right|_g \\ \left|\nabla^g_X \mathscr{A}^{-1/2}\right|_g &\leqslant |\mathscr{A}|\,\left|\mathscr{A}^{-1/2}\right|^3\,\left|\nabla^h_X - \nabla^g_X\right|_g , \end{aligned}$$

*where $|\cdot|_g$ denotes the operator norm induced by the inner product $g$.*

**Proof.** We divide the proof into three steps.

**1°** Recall that $\mathscr{A}^{1/2}$ is an isometry, so that

$$\left(\mathscr{A}^{1/2}\alpha, \mathscr{A}^{1/2}\alpha\right)_g = (\alpha, \alpha)_h\,.$$

If we differentiate the identity in the direction of $X$, we get on the one hand

$$X\left(\mathscr{A}^{1/2}\alpha, \mathscr{A}^{1/2}\alpha\right)_g = 2\left(\nabla^g_X\mathscr{A}^{1/2}(\alpha), \mathscr{A}^{1/2}\alpha\right)_g + 2\left(\mathscr{A}^{1/2}\nabla^g_X\alpha, \mathscr{A}^{1/2}\alpha\right)_g,$$

whereas, we have on the other hand

$$X(\alpha, \alpha)_h = 2\left(\nabla^h_X\alpha, \alpha\right)_h = 2\left(\mathscr{A}\nabla^h_X\alpha, \alpha\right)_g.$$

Thus,

$$\left(\nabla^g_X\mathscr{A}^{1/2}\alpha, \mathscr{A}^{1/2}\alpha\right)_g = \left(\mathscr{A}^{1/2}\circ\left(\nabla^h_X - \nabla^g_X\right)\alpha, \mathscr{A}^{1/2}\alpha\right)_g. \tag{2.7}$$

**2°** Clearly, as $\mathscr{A}^{1/2}$ is self-adjoint so is $\nabla^g_X\mathscr{A}^{1/2}$. So, let $\lambda$ be an eigenvalue of $\nabla^g_X\mathscr{A}^{1/2}$ with $|\lambda| = \left|\nabla^g_X\mathscr{A}^{1/2}\right|_g$ at a fixed point $x \in M$. Let $v \in \mathsf{T}_xM$ be a corresponding $g$-normalised eigenvector. We estimate

$$\left|\left(\mathscr{A}^{1/2}\nabla^g_X\mathscr{A}^{1/2}v, v\right)_g\right| = |\lambda|\left|\left(\mathscr{A}^{1/2}v, v\right)_g\right| = \left|\nabla^g_X\mathscr{A}^{1/2}\right|_g\left|\left(\mathscr{A}^{1/2}v, v\right)_g\right| \geqslant \left|\nabla^g_X\mathscr{A}^{1/2}\right|_g\left|\mathscr{A}^{-1/2}\right|^{-1},$$

where we used that $\mathscr{A}$, hence $\mathscr{A}^{1/2}$, has only (strictly) positive eigenvalues: More precisely, since $\mathscr{A}^{1/2}$ is self-adjoint, there are eigenvalues $\lambda_{\min}$ and $\lambda_{\max}$ of $\mathscr{A}^{1/2}$ such that

$$\lambda_{\min}(\mathscr{A}^{1/2}) \leqslant \left(\mathscr{A}^{1/2}u, u\right)_g \leqslant \lambda_{\max}(\mathscr{A}^{1/2}).$$

By definition, the smooth vector bundle morphism $\mathscr{A}$ has only strictly positive eigenvalues, so does $\mathscr{A}^{1/2}$ and we have that $\left(\mathscr{A}^{1/2}u, u\right)_g > 0$. Moreover,

$$\left|\mathscr{A}^{1/2}\right| \geqslant \lambda_{\min}(\mathscr{A}^{1/2}) = \frac{1}{\lambda_{\max}(\mathscr{A}^{-1/2})} = \left|\mathscr{A}^{-1/2}\right|^{-1}.$$

Hence the first estimate follows:

$$\begin{aligned}\left|\nabla^g_X\mathscr{A}^{1/2}\right|_g &\leqslant \left|\mathscr{A}^{-1}\right|^{1/2}\sup_{|v|\leqslant 1}\left|\left(\mathscr{A}^{1/2}\nabla^g_X\mathscr{A}^{1/2}v, v\right)_g\right| \\ &\overset{(2.7)}{=} \left|\mathscr{A}^{-1}\right|^{1/2}\sup_{|v|\leqslant 1}\left|\left(\mathscr{A}\circ\left(\nabla^h_X - \nabla^g_X\right)v, v\right)_g\right| \\ &\leqslant \left|\mathscr{A}^{-1}\right|^{1/2}|\mathscr{A}|\left|\nabla^h_X - \nabla^g_X\right|_g.\end{aligned}$$

**3°** Covariantly differentiating the identity

$$\mathrm{id} = \mathscr{A}^{1/2}\circ\mathscr{A}^{-1/2},$$

we get

$$0 = \nabla^g_X\mathscr{A}^{1/2}\circ\mathscr{A}^{-1/2} + \mathscr{A}^{1/2}\nabla^g_X\mathscr{A}^{-1/2}.$$

Thus,

$$\nabla^g_X\mathscr{A}^{-1/2} = -\mathscr{A}^{-1/2}\nabla^g_X\mathscr{A}^{1/2}\circ\mathscr{A}^{-1/2},$$

and the second estimate follows from by part **2°**. ∎

Denoting by $0 < \rho_{g,h} =: \rho \in \mathrm{C}^\infty(M)$ the Radon-Nikodým density, i.e.

$$\mathrm{d}\,\mathrm{vol}_h = \rho_{g,h} \mathrm{d}\,\mathrm{vol}_g,$$

by quasi-isometry, the following identities hold

$$\rho_{h,g} = \rho_{g,h}^{-1}, \quad \mathcal{A}_{h,g} = \mathcal{A}_{g,h}^{-1}, \quad \rho_{g,h} = (\det A)^{1/2}, \quad 0 < \inf \rho_{g,h} \leqslant \sup \rho_{g,h} < \infty. \tag{2.8}$$

We now define the zeroth order deviation of the two metrics (considered as multiplicative perturbations of each other) as

$$\delta_{g,h}(x) := 2 \sinh\left( \frac{m}{4} \max_{\lambda \in \sigma(A_{g,h}(x))} |\log \lambda| \right) = \max_{\lambda \in \sigma(A_{g,h}(x))} \left| \lambda^{\frac{m}{4}} - \lambda^{-\frac{m}{4}} \right| : M \to (0, \infty), \tag{2.9}$$

symmetric in $g$ and $h$ by quasi-isometry, i.e. $\delta_{g,h} = \delta_{h,g}$. Note that these deviations have been considered first in [HPW14, Appendix A]. We will make use of the fact that

$$\sup \delta_{g,h}(x) < \infty \quad \Longleftrightarrow \quad g \sim h.$$

The definition is becoming clearer in the Proof of Lemma 5.2. Moreover, let

$$\delta^{\nabla}_{g,h}(x) := \left| \nabla^h - \nabla^g \right|_g^2 (x) : M \to [0, \infty) \tag{2.10}$$

be a weight function defined in terms of the corresponding covariant derivatives $\nabla^h$ and $\nabla^g$, defined in terms of the operator norm induced by the inner product $g$.

**Remark 2.6.** Note that the difference of two connections $\nabla^h - \nabla^g$ is a one-form on $M$ with values in $\mathrm{End}\,\mathrm{T}M$, i.e.

$$\nabla^h - \nabla^g \in \Gamma(\mathrm{T}^* M \otimes \mathrm{End}\,\mathrm{T}M).$$

**Example 2.7** (Conformal metric change). Let $h := g_\psi$ be a conformal perturbation of $g$, i.e. we set $g_\psi := \mathrm{e}^{2\psi} g$ for some smooth function $\psi : M \to \mathbb{R}$. We take $A := e^{2\psi}$, so $A^{-1} = \mathrm{e}^{-2\psi}$ and

$$\delta_{g,h}(x) = 2 \sinh \frac{m}{4} |\psi(x)|.$$

Hence,

$$g \sim h \quad \Longleftrightarrow \quad \psi \quad \text{bounded.}$$

By (6.2d), for any smooth vector field $X, Y \in \Gamma_{\mathrm{C}^\infty}(\mathrm{T}M)$, we have

$$\left( \nabla^{g_\psi}_X - \nabla^g_X \right) Y = \mathbf{d}\psi(X) Y + \mathbf{d}\psi(Y) X - (X, Y)_g \operatorname{grad}_g \psi. \tag{2.11}$$

By norm equivalence on finite dimensional spaces, let us work with the Hilbert-Schmidt norm for the calculation. Let $(X_i)_{i=1}^m$ a smooth local $g$-orthonormal frame of vector fields. Then in local coordinates,

$$\begin{aligned}|\mathbf{d}\psi|^2_{g\otimes \mathrm{HS}} &= \sum_{i=1}^m \left|X_i\psi\right|^2,\\ \left|\nabla^{g_\psi}-\nabla^g\right|^2_{g\otimes\mathrm{HS}} &= \sum_{j,k=1}^m \left|\left(\nabla^{g_\psi}-\nabla^g\right)(X_j,X_k)\right|^2\\ &= \sum_{j,k=1}^m \left|(X_j\psi)X_k+(X_k\psi)X_j-\delta_{jk}\sum_i (X_i\psi)X_i\right|^2\\ &= \sum_{j<k} 2\left|(X_j\psi)X_k+(X_k\psi)X_j\right|^2+\sum_{i=1}^m\left|X_i\psi\right|^2,\end{aligned}$$

so that

$$|\mathbf{d}\psi|^2_g \lesssim \left|\nabla^{g_\psi}-\nabla^g\right|^2_g = \delta^\nabla_{g,g_\psi}.$$

Next, we give a formula how to express the codifferential $\boldsymbol{\delta}_h$ with respect the metric $h$ in terms of the codifferential $\boldsymbol{\delta}_g$ in terms of $g$.

**Lemma 2.8.** *Let $g,h\in \mathrm{Metr}M$, $g\sim h$, and $\mathscr{A} := \mathscr{A}_{g,h}$ be the smooth vector bundle morphism defined by (2.4). Then the codifferential with respect to the metric $h$ is given by*

$$\boldsymbol{\delta}_h\eta = \mathscr{A}^{-1}\left(\boldsymbol{\delta}_g(\mathscr{A}\eta)-(\mathbf{d}\log\rho)^{\sharp_g}\lrcorner_g(\mathscr{A}\eta)\right)\qquad\text{for all }\ \eta\in\Omega_{C^\infty}(M).$$

**Proof.** For any $\eta_1\in\Omega_{C^\infty}(M)$, $\eta_2\in\Omega_{C^\infty}(M)$, we calculate

$$\begin{aligned}\left\langle\eta_1,\mathscr{A}^{-1}\left(\boldsymbol{\delta}_g(\mathscr{A}\eta_2)-(\mathbf{d}\log\rho)^{\sharp_g}\lrcorner_g(\mathscr{A}\eta_2)\right)\right\rangle_h &= \left\langle\rho\mathscr{A}\eta_1,\mathscr{A}^{-1}\left(\boldsymbol{\delta}_g(\mathscr{A}\eta_2)-\frac{(\mathbf{d}\rho)^{\sharp_g}}{\rho}\lrcorner_g(\mathscr{A}\eta_2)\right)\right\rangle_g\\ &= \left\langle\eta_1,\rho\boldsymbol{\delta}_g(\mathscr{A}\eta_2)-(\mathbf{d}\rho)^{\sharp_g}\lrcorner_g(\mathscr{A}\eta_2)\right\rangle_g\\ &= \langle\mathbf{d}(\rho\eta_1),\mathscr{A}\eta_2\rangle_g-\langle\mathbf{d}\rho\wedge\eta_1,\mathscr{A}\eta_2\rangle_g\\ &= \langle\mathbf{d}\rho\wedge\eta_1,\mathscr{A}\eta_2\rangle_g+\langle\rho\mathbf{d}\eta_1,\mathscr{A}\eta_2\rangle_g-\langle\mathbf{d}\rho\wedge\eta_1,\mathscr{A}\eta_2\rangle_g\\ &= \langle\mathbf{d}\eta_1,\eta_2\rangle_h,\end{aligned}$$

where we used that $\mathscr{A}$ is fibrewise self-adjoint. ■

In the Proof of the main result, the gradient of the logarithm of the Radon-Nikodým density $\rho_{g,h}$ can be estimated in terms of smooth vector bundle morphism $\mathscr{A}_{g,h}$ and $\delta^\nabla_{g,h}$ which is reflected in the next Proposition 2.12. Therefore, we note two auxiliary lemmas.

**Lemma 2.9.** *Let $A,B$ be two complex $m\times m$-matrices. Then,*

$$|\mathrm{tr}(AB)|\leqslant \|A\|_{\mathrm{HS}}\,\|B\|_{\mathrm{HS}},$$

*where $\|\cdot\|_{\mathrm{HS}}$ denotes the Hilbert-Schmidt norm.*

**Proof.** Let $A, B$ be two complex $m \times m$-matrices with singular values

$$\sigma_1(A) \geqslant \sigma_2(A) \geqslant \ldots \geqslant \sigma_m(A) \qquad \text{of } A \text{ and}$$
$$\sigma_1(B) \geqslant \sigma_2(B) \geqslant \ldots \geqslant \sigma_m(B) \qquad \text{of } B.$$

Then, with the Hilbert-Schmidt norm $\|\cdot\|_{\mathrm{HS}}$,

$$\|\sigma(A)\| = \sqrt{\sum_{i=1}^{m} \sigma_i^2(A)} = \sqrt{\operatorname{tr}(A^* A)} =: \|A\|_{\mathrm{HS}} \tag{2.12}$$

by norm equivalence on finite-dimensional vector spaces. By the well-known von Neumann trace formula [Mir75], we get

$$\begin{aligned} |\operatorname{tr}(AB)| &\leqslant \sum_{i=1}^{m} \sigma_i(A)\sigma_i(B) \\ &= (\sigma(A), \sigma(B))_{\mathbb{R}^m} \\ &\leqslant \|\sigma(A)\| \;\; \|\sigma(B)\| = \|A\|_{\mathrm{HS}} \, \|B\|_{\mathrm{HS}} \, , \end{aligned}$$

where we used Cauchy-Schwarz for the first inequality and (2.12) in the last equality. ■

**Lemma 2.10** (Classical Jacobi's formula). *Let $B = B(t)$ an $m \times m$-matrix parametrised by $t$. If $B$ is invertible, then*

$$\frac{\mathrm{d}}{\mathrm{d}t} \det B(t) = \det B(t) \ \operatorname{tr}\left( B(t)^{-1} \frac{\mathrm{d}}{\mathrm{d}t} B(t) \right) . \tag{2.13}$$

Next, we extend the classical Jacobi formula to our setting.

**Lemma 2.11** (Jacobi's formula). *Let $A := A_{g,h}$ be the smooth vector bundle morphism defined by (2.3) ($x \in M$). Then, for any $X \in \Gamma(\mathrm{T}M)$,*

$$\mathbf{d}(\det A)X = \det A \ \operatorname{tr}\left(A^{-1} \nabla_X A\right) .$$

**Proof.** For $\varepsilon > 0$, let $\gamma : (-\varepsilon, \varepsilon) \to M$ be a smooth curve on $M$ such that $\gamma(0) = x$ and $\dot\gamma(0) = X$. Then, computing the differential using a velocity vector,

$$\begin{aligned} \mathbf{d}(\det A)_x X = \frac{\mathrm{d}}{\mathrm{d}t} \det A(\gamma(t))\Big|_{t=0} &= \det A(\gamma(t)) \ \operatorname{tr}\left( A(\gamma(t))^{-1} \frac{\mathrm{d}}{\mathrm{d}t} A(\gamma(t)) \right)\Big|_{t=0} \\ &= \det A(\gamma(t)) \ \operatorname{tr}\left( A(\gamma(t))^{-1} \, \nabla_{\dot\gamma(t)} A(\gamma(t)) \right)\Big|_{t=0} \\ &= \det A(x) \ \operatorname{tr}\left( A(x)^{-1} \nabla_X A(x) \right) , \end{aligned}$$

using Lemma 2.10 in the second step. ■

We are now in the position to prove the aimed estimate of the the Radon-Nikodým density.

**Proposition 2.12.** *Let $g, h \in \mathrm{Metr}M$, $g \sim h$, and $\mathscr{A} := \mathscr{A}_{g,h}$ be the smooth vector bundle morphism defined by (2.4) and $\rho = (\det A)^{1/2}$ the Radon-Nikodým density defined by (2.8). Then we can estimate as follows:*

$$|\mathbf{d} \log \rho|_g \leqslant C(m) \left| \mathscr{A}^{-1} \right| |\nabla^g \mathscr{A}|_g \, .$$

**Proof.** We first recall, that for $x \in M$,

$$\dim \bigwedge^k \mathsf{T}_x^* M = \binom{m}{k} \qquad \Longrightarrow \qquad \dim \bigwedge \mathsf{T}_x^* M = 2^m < \infty,$$

so $\dim \bigwedge \mathsf{T}_x^* M$ is finite-dimensional ($x \in M$). Recall that the Radon-Nikodým density $\rho$ is given by $\rho = (\det A)^{1/2}$. By Jacobi's formula, Lemma 2.11 above,

$$\mathbf{d}(\log \rho)_x X = \frac{1}{2} \operatorname{tr}\left(A(x)^{-1} \nabla_X A(x)\right).$$

Hence, using Lemma 2.9,

$$\begin{aligned} \left| \mathbf{d}(\log \rho)_x X \right| &\leqslant \frac{1}{2} \left\| A(x)^{-1} \right\|_{\mathrm{HS}} \left\| \nabla_X A(x) \right\|_{g \otimes \mathrm{HS}} \\ &\leqslant C(m) \left| A(x)^{-1} \right| \left| \nabla_X A(x) \right|_g \\ &\leqslant C(m) \left| \mathscr{A}(x)^{-1} \right| \left| \nabla_X \mathscr{A}(x) \right|_g \end{aligned}$$

by norm equivalence on finite-dimensional spaces. ■

Next, we define the bounded identification operator

$$\begin{aligned} I = I_{g,h} : \Omega_{\mathsf{L}^2}(M, g) &\to \Omega_{\mathsf{L}^2}(M, h) \\ \alpha(x) &\mapsto \mathscr{A}_{g,h}^{-1/2}(x) \alpha(x), \end{aligned}$$

well-defined by $g \sim h$.

**Lemma 2.13.** *Let $g, h \in \mathsf{Metr} M$, $g \sim h$. The adjoint $I^*$ of the bounded identification operator $I$ is given by*

$$\begin{aligned} I^* = I_{g,h}^* : \Omega_{\mathsf{L}^2}(M, h) &\to \Omega_{\mathsf{L}^2}(M, g) \\ \alpha(x) &\mapsto \rho_{g,h}(x) \mathscr{A}_{g,h}^{1/2}(x) \alpha(x). \end{aligned} \tag{2.14}$$

**Proof.** For compactly supported $\alpha \in \Omega(M, h)$ and $\beta \in \Omega(M, g)$, we get

$$\begin{aligned} \left\langle I_{g,h}^* \alpha, \beta \right\rangle_{\Omega_{\mathsf{L}^2}(M,g)} &= \left\langle \alpha, I_{g,h} \beta \right\rangle_{\Omega_{\mathsf{L}^2}(M,h)} \\ &= \int_M \left( \alpha, \mathscr{A}_{g,h}^{-1/2} \beta \right)_h \mathrm{d}\,\mathrm{vol}_h \\ &= \int_M \left( \mathscr{A}_{g,h}^{1/2} \alpha, \mathscr{A}_{g,h}^{1/2} \mathscr{A}_{g,h}^{-1/2} \beta \right)_g \rho_{g,h} \mathrm{d}\,\mathrm{vol}_g \\ &= \int_M \left( \rho_{g,h} \mathscr{A}_{g,h}^{1/2} \alpha, \beta \right)_g \mathrm{d}\,\mathrm{vol}_g = \left\langle \rho_{g,h} I_{g,h}^{-1} \alpha, \beta \right\rangle_{\Omega_{\mathsf{L}^2}(M,g)}. \end{aligned}$$

■

For every $g \in \mathsf{Metr} M$,

$$(P_s^g)_{s>0} := \left( \mathrm{e}^{\frac{s}{2} \Delta_g} \right)_{s>0} \subset \mathscr{L}(\Omega_{\mathsf{L}^2}(M, g))$$

is the heat semigroup defined by the spectral theorem, choosing $f : \mathbb{R} \to \mathbb{R}$, $f(\lambda) := \mathrm{e}^{\lambda t/2}$. Let us denote by

$$(0, \infty) \times M \times M \ni (s, x, y) \mapsto p_s^g(x, y) := \mathrm{e}^{\frac{s}{2} \Delta_g}(x, y) \in \mathsf{Hom}\left( \bigwedge \mathsf{T}_y^* M, \bigwedge \mathsf{T}_x^* M \right)$$

the corresponding jointly smooth integral kernel of $P_s^g$. The smooth representative of $P_s\alpha(x)$ is given by

$$P_s^g\alpha(x) = \int_M e^{\frac{s}{2}\Delta_g}(x,y)\alpha(y)\,\mathrm{vol}_g(\mathrm{d}y).$$

It is well-known [Gün17, Theorem II.1.] that we have

$$\int_M \left|p_s^g(x,y)\right|_g^2 \mathrm{vol}_g(\mathrm{d}y) < \infty, \qquad \text{for all } \; s > 0, \; x \in M.$$

The form degree $k$ will be indicated again by round brackets $p_s^{g,(k)}$. Finally, we set

$$\Phi_g(x,s) := \sup_M p_s^{g,(0)}(x,\cdot), \tag{2.15}$$

indicating the minimal heat kernel $p_s^{g,(0)}$ acting on 0-forms, i.e. functions. Then it is well-known that, for all $(x,s) \in M \times (0,\infty)$, it follows that $\Phi_g(x,s) < \infty$. One can even show [Gün17] that

$$\sup_K \Phi_g(\cdot,s) < \infty \qquad \text{for all } \; s > 0, \; K \subset M \text{ compact.}$$

## 3 Bismut-type Formulae & Gradient Estimates

In this section, we will *omit the metric $g$ in the notation if unambiguous*.

For a brief and concise overview to basic notions in stochastic differential geometry and stochastic analysis on manifolds, we would like to refer the reader to e.g. [Tha16] and [Hsu02]. Being a German reader, a profound introduction can be also found in [HT94].

Let $(\Omega, \mathscr{F}, (\mathscr{F}_t)_{t\geqslant 0}, \mathbb{P})$ be a filtered probability space satisfying the usual hypotheses, i.e. a probability space $(\Omega, \mathscr{F}, \mathbb{P})$ equipped with a filtration that is right-continuous

$$\mathscr{F}_t = \mathscr{F}_{t+} := \bigcap_{s>t} \mathscr{F}_s \qquad \text{for all } t \geqslant 0,$$

and complete, i.e. $\mathscr{F}_0$ contains all subsets of $\mathbb{P}$-null sets.

Let $X(x)$ be a Brownian motion on $M$ starting at $x \in M$ and $\zeta(x)$ its maximal lifetime, i.e. an $M$-valued stochastic process $\left(X_t = X_t(x)\right)_{0\leqslant t<\zeta}$ such that $X_0 = x \in M$ and

$$N_t^f := f(X_t) - f(X_0) - \frac{1}{2}\int_0^t \Delta_M f(X_s)\,\mathrm{d}s \quad \text{ on } \{t < \zeta(x)\}$$

is a real local martingale for all $f \in C^\infty(M)$ and $\Delta_M$ denoting the Laplace-Beltrami operator on $M$. In other words, the process $X(x)$ is a diffusion process in $M$ starting at $x \in M$ which is generated by $\frac{1}{2}\Delta_M$.

Further, let $/\!/$ be the stochastic parallel transport along $X(x)$ and $B$ the stochastic anti-development of $X$ to $\mathrm{T}_xM$, which is a standard Brownian motion on $\mathrm{T}_xM \cong \mathbb{R}^m$. For every $r > 0$ let

$$\tau(x,r) := \inf\left\{t \geqslant 0 \,:\, X_t(x) \notin \mathrm{B}(x,r)\right\} \,:\, \Omega \to [0,\infty]$$

be the first exit time of $X$ from the open ball $\mathrm{B}(x,r)$.

The expected value of $X$ will be denoted by $\mathbb{E}X := \mathbb{E}(X)$. To relax notation, we simply write $\mathbb{E}f(X)$ rather than $\mathbb{E}(f(X))$ and $\left(\mathbb{E}X^2\right)^{1/2}$ rather than $\left(\mathbb{E}(X^2)\right)^{1/2}$ etc.

Let $E, \widetilde{E}$ be two Riemannian vector bundles over $M$, endowed with a metric connection $\nabla^E$ and $\nabla^{\widetilde{E}}$ respectively. For a given multiplication map $m \in \Gamma(\mathrm{Hom}(\mathrm{T}^*M \otimes E, \widetilde{E})) \cong \Gamma(\mathrm{T}M \otimes E^* \otimes \widetilde{E})$, we consider the Dirac-type operator

$$\mathsf{D}_m := m\nabla^E : \Gamma(E) \to \Gamma(\widetilde{E})$$

which is understood as the composition

$$\Gamma(E) \xrightarrow{\nabla^E} \Gamma(\mathrm{T}^*M \otimes E) \xrightarrow{m} \Gamma(\widetilde{E}).$$

A multiplication map $m$ is said to be compatible with $\nabla$ provided $\nabla m = 0$, i.e.

$$\nabla_v^{\widetilde{E}}(m_U \alpha) = m_{\nabla_v^{\mathrm{T}M} U}\alpha + m_U(\nabla_v^E \alpha) \qquad \text{for all } U \in \Gamma(\mathrm{T}M),\ \alpha \in \Gamma(E),\ v \in \mathrm{T}M.$$

The horizontal Laplacian $\square$ is the second order differential operator given by the following decomposition

$$\square : \Gamma(E) \xrightarrow{\nabla^E} \Gamma(\mathrm{T}^*M \otimes E) \xrightarrow{\nabla^{\mathrm{T}^*M \otimes E}} \Gamma(\mathrm{T}^*M \otimes \mathrm{T}^*M \otimes E) \xrightarrow{\mathrm{tr}} \Gamma(E).$$

Driver and Thalmaier [DT01, p. 48] propose the following formalism: Let $\mathsf{L}$ and $\widetilde{\mathsf{L}}$ are given second order differential operators on $\Gamma(E)$ and $\Gamma(\widetilde{E})$ respectively that satisfy the following two conditions.

(1) The operators $\mathsf{D}_m$, $\mathsf{L}$ and $\widetilde{\mathsf{L}}$ obey the commutation rule, for some $\varrho \in \Gamma(\mathrm{Hom}(E, \widetilde{E}))$,

$$\mathsf{D}_m \mathsf{L} = \widetilde{\mathsf{L}}\mathsf{D}_m - \varrho.$$

(2) The operators $\mathscr{R} := \square - \mathsf{L} : \Gamma(E) \to \Gamma(E)$ and $\widetilde{\mathscr{R}} := \widetilde{\square} - \widetilde{\mathsf{L}} : \Gamma(\widetilde{E}) \to \Gamma(\widetilde{E})$ are zeroth order operators, i.e. $\mathscr{R}$ and $\widetilde{\mathscr{R}}$ are section in $\Gamma(\mathrm{End}\, E)$ and $\Gamma(\mathrm{End}\, \widetilde{E})$, provided $m$ is compatible with the Levi-Civita connection.

In geometrically natural situations, we have $\varrho = 0$ or $\rho \in \Gamma(\mathrm{Hom}(E, \widetilde{E}))$ is of zeroth order. Under those assumptions, Driver and Thalmaier [DT01] can prove derivative formulae for the heat semigroup in the general setting of vector bundles using martingale methods. For a detailed discussion we refer the reader to [DT01]. Let us note two important examples.

**Example 3.1.** The exterior bundle of total forms $E = \bigwedge \mathrm{T}^*M \to M$ with its natural connection

$$\nabla := \nabla^{\bigwedge \mathrm{T}^*M} := \bigoplus_{k=0}^{m} \nabla^{\bigwedge^k \mathrm{T}^*M}$$

and Clifford action $c : \mathrm{T}M \to \mathrm{End}(\bigwedge \mathrm{T}^*M)$, $c(\alpha)\beta := \alpha \wedge \beta - \alpha^\sharp \lrcorner \beta$. The Hodge Laplacian $\Delta$ is related to the horizontal Laplacian $\square$ by the Weitzenböck formula

$$\Delta = \square - \mathscr{R}, \tag{3.1}$$

where Weitzenböck curvature operator $\mathscr{R} \in \Gamma(\operatorname{End} \Omega_{C^\infty}(M))$ is a symmetric field of endomorphisms. Acting on $k$-forms, the field of endomorphisms is specified again by an index

$$\mathscr{R}^{(k)} := \mathscr{R}\big|_{\Omega^k_{C^\infty}(M,g)}.$$

In particular, note that $\mathscr{R}^{(1),\mathrm{tr}} = \mathrm{Ric}$ and $\mathscr{R}^{(0),\mathrm{tr}} = 0$. Moreover, it can be written explicitly (cf. e.g. [DT01, Lemma A.7]), for any orthonormal basis $(e_k)_{1\leqslant k\leqslant m}$,

$$\mathscr{R}^{(k)} = \sum_{k,l=1}^{m} \mathrm{R}(e_l, e_k)(e_l^\flat \wedge \bullet)(e_k \lrcorner \bullet),$$

where $\mathrm{R}(e_l, e_k)$ is the curvature tensor acting on $k$-forms (cf. [DT01, Lemma A.9]). Then

$$\mathbf{D}_g \equiv \mathbf{D}_c = \mathbf{d} + \boldsymbol{\delta}, \qquad \mathbf{L} = \widetilde{\mathbf{L}} = \Delta, \qquad m = c \quad \text{and} \quad \varrho = 0.$$

In particular, for $E := \bigwedge^k \mathrm{T}^*M$ and $\widetilde{E} := \bigwedge^{k+1} \mathrm{T}^*M$, then $\quad \varrho = 0 \quad$ with

$$\mathbf{D}_m = \mathbf{d}\big|_{\Omega^k}, \qquad \mathbf{L} = -\Delta^{(k)}, \quad \widetilde{\mathbf{L}} = \Delta^{(k+1)}, \qquad \mathscr{R} = \mathscr{R}^{(k)}, \quad \widetilde{\mathscr{R}} = \mathscr{R}^{(k+1)}, \qquad m(\alpha \otimes \beta) = \alpha \wedge \beta.$$

If instead $\widetilde{E} := \bigwedge^{k-1} \mathrm{T}^*M$, then again $\quad \varrho = 0 \quad$ but with

$$\mathbf{D}_m = \boldsymbol{\delta}\big|_{\Omega^k}, \qquad \mathbf{L} = \Delta^{(k)}, \quad \widetilde{\mathbf{L}} = \Delta^{(k-1)}, \qquad \mathscr{R} = \mathscr{R}^{(k)}, \quad \widetilde{\mathscr{R}} = \mathscr{R}^{(k-1)}, \qquad m(\alpha \otimes \beta) = -(\alpha^\sharp \lrcorner \beta).$$

**Example 3.2** ([DT01, cf. Proposition 2.15]). Let $\widetilde{E} = \mathrm{T}^*M \otimes E$ and $m = \mathrm{id}_{\widetilde{E}}$. For $\mathbf{D}_m = \nabla$ and given $\mathscr{R} \in \operatorname{End} E$, we set

$$\widetilde{\mathscr{R}} = \mathrm{Ric}^{\mathrm{tr}} \otimes \mathbf{1}_E - 2\mathrm{R}^E \cdot + \mathbf{1}_{\mathrm{T}^*M} \otimes \mathscr{R} \in \Gamma(\operatorname{End} \widetilde{E}), \tag{3.2}$$

$$\varrho = \nabla \cdot \mathrm{R}^E + \nabla^{\operatorname{End} E} \mathscr{R} \in \Gamma(\operatorname{Hom}(E, \widetilde{E})), \tag{3.3}$$

where $\mathrm{R}^E$ denotes the Riemannian curvature tensor to $\nabla$ on $E$ and

$$\mathrm{Ric}^{\mathrm{tr}} \in \Gamma(\operatorname{End} \mathrm{T}^*M)$$

denotes the transpose of the Ricci curvature tensor $\mathrm{Ric} \in \Gamma(\operatorname{End} \mathrm{T}M)$ on $M$, and for any $\eta \in \widetilde{E}_x$, $v \in \mathrm{T}_xM$, $\alpha \in E_x$, $(e_i)$ an orthonormal frame for $\mathrm{T}_xM$,

$$\left(\mathrm{R}^E \cdot \eta\right)(v) := \sum_{i=1}^{m} \mathrm{R}^E(v, e_i)\eta(e_i),$$

$$\left(\nabla \cdot \mathrm{R}^E \alpha\right)(v) := \sum_{i=1}^{m} \nabla_{e_i} \mathrm{R}^E(e_i, v)\alpha,$$

$$(\nabla \mathscr{R} \alpha)(v) := (\nabla_v \mathscr{R})\alpha.$$

Choosing $m = \mathrm{id}$, for $\mathbf{L} := \square - \mathscr{R}$, $\widetilde{\mathbf{L}} = \widetilde{\square} - \widetilde{\mathscr{R}}$, it follows that $\varrho \in \Gamma(\operatorname{Hom}(E, \widetilde{E}))$.

Next, we introduce the stochastic representation of the semigroup. To this end, recall that $/\!/_t^E : E_x \to E_{X_t}$ is the parallel transport along our diffusion $X_t = X_t(x)$ started at $x \in M$. We define the linear operators on $E_x$ and $\widetilde{E}_x$, respectively,

$$\mathscr{R}_{/\!/_t} := (/\!/_t^E)^{-1} \mathscr{R} \, /\!/_t^E \qquad \text{and} \qquad \widetilde{\mathscr{R}}_{/\!/_t} := (/\!/_t^{\widetilde{E}})^{-1} \widetilde{\mathscr{R}} \, /\!/_t^{\widetilde{E}}$$

along the paths of $X(x)$ in the following way: Via the stochastic parallel transport we get to a random point on the tangent space at $X_t$ and apply the curvature (in case $k = 1$ just Ric and $E = \mathrm{T}M$) considered as a linear transformation. Then we parallel transport back to where the diffusion started. Thus, we get a linear mapping $E_x \to E_x$ which now depends on random, i.e.

$$\begin{array}{ccc} E_x & \overset{\mathscr{R}_{/\!/_t^E}}{\dashrightarrow} & E_x \\ {\scriptstyle /\!/_t^E}\downarrow & & \uparrow{\scriptstyle /\!/_t^{E,-1}} \\ E_{X_t} & \underset{\mathscr{R}_{X_t}}{\longrightarrow} & E_{X_t} \end{array}$$

Let $\mathscr{Q}_t$ be the End$(E_x)$-valued, and $\widetilde{\mathscr{Q}}_t$ the End$(\widetilde{E}_x)$-valued respectively, pathwise solutions to the ordinary differential equations

$$\begin{aligned} \frac{\mathrm{d}}{\mathrm{d}t}\mathscr{Q}_t &= -\frac{1}{2}\mathscr{R}_{/\!/_t}\mathscr{Q}_t, \qquad \mathscr{Q}_0 = \mathrm{id}_{E_x}, \\ \frac{\mathrm{d}}{\mathrm{d}t}\widetilde{\mathscr{Q}}_t &= -\frac{1}{2}\widetilde{\mathscr{R}}_{/\!/_t}\widetilde{\mathscr{Q}}_t, \qquad \widetilde{\mathscr{Q}}_0 = \mathrm{id}_{\widetilde{E}_x}. \end{aligned} \tag{3.4}$$

The composition $\mathscr{Q} \circ /\!/^{E,-1}$ maps from a random point $X_t$ back to the starting point $x$ and is called the (inverse) **damped parallel transport** along the paths of $X(x)$.

Let $\mathscr{H}_1$ and $\mathscr{H}_2$ be two Hilbert spaces. By $\mathscr{L}(\mathscr{H}_1, \mathscr{H}_2)$ we denote the space of all bounded operators $P : \mathscr{H}_1 \to \mathscr{H}_2$. For any $q \in [1, \infty)$, the operator $P$ is a **Schatten operator of class** $q$ if $\operatorname{tr}\left|P^*P\right|^q < \infty$. For $q = 1$ the operator $P$ is trace class and for $q = 2$ the operator $P$ is in the Hilbert-Schmidt class. Given two metric vector bundles $E$, $\widetilde{E}$ over $M$ and a bounded operator

$$P \in \mathscr{L}(\Gamma_{\mathrm{L}^2}(M, E), \Gamma_{\mathrm{L}^2}(M, \widetilde{E}))$$

such that there exists a corresponding jointly smooth integral kernel $p(x, y)$ of $P$

$$M \times M \ni (x, y) \mapsto p_s(x, y) \in \mathrm{Hom}\left(E_y, \widetilde{E}_x\right),$$

the uniquely determined map such that we have

$$P\alpha(x) = \int_M p(x, y)\alpha(y)\,\mathrm{vol}_\bullet(\mathrm{d}y).$$

Then $P$ is Hilbert-Schmidt, if

$$\iint_{M\times M} |p(x, y)|^2\,\mathrm{vol}_\bullet(\mathrm{d}x)\,\mathrm{vol}_\bullet(\mathrm{d}y) < \infty.$$

In particular, it is well-known [RS79] that Hilbert-Schmidt operators and the product of a bounded operator with some Hilbert-Schmidt operator are compact, and that the product of two Hilbert-Schmidt class operators is trace class.

Given a potential $w : M \to \mathbb{C}$, the Feynman-Kac semigroup

$$P_s^w f(x) := \mathbb{E}\left(\mathrm{e}^{-\frac{1}{2}\int_0^s w(X_r(x))\mathrm{d}r} f(X_s(x))\mathbb{1}_{\{s<\zeta(x)\}}\right)$$

acts on (bounded) measurable functions $f$ on $M$. Further, let $\underline{\mathscr{R}}^E = \sigma_{\min}(\mathscr{R}^E)$, i.e.

$$\underline{\mathscr{R}}^E(x) := \min\left\{\left(\mathscr{R}_x v, v\right), v \in E_x, |v| = 1\right\}.$$

By uniform continuity, $\underline{\mathscr{R}}^E$ is a continuous function on $M$. By Gronwall's inequality,

$$\left|\mathcal{Q}_s\right|_{\mathrm{op}} \leqslant \exp\left(-\frac{1}{2}\int_0^s \underline{\mathscr{R}}^E(X_r(x))\mathrm{d}r\right),$$

where $|\cdot|_{\mathrm{op}}$ is the usual operator norm. Then, we have the following probabilistic representation of the semigroup, for all $s > 0$ and every $x \in M$,

$$P_s^g \alpha(x) = \int_M p_s^g(x, y)\alpha(y)\,\mathrm{vol}_g(\mathrm{d}y) = \mathbb{E}\left(\mathcal{Q}_s /\!/_s^{E,-1} \alpha(X_s(x))\mathbb{1}_{\{s<\zeta(x)\}}\right) \qquad \text{for all} \quad \alpha \in \Gamma_{L^2}(E),$$

provided the scalar semigroup $P_t^{\underline{\mathscr{R}}^E} |\alpha|(x) < \infty$. In particular, we get semigroup domination

$$\left|P_s^g \alpha(x)\right| \leqslant P_s^{\underline{\mathscr{R}}^E} |\alpha|(x).$$

The existence of the Feynman-Kac semigroup will be essential to prove the Bismut-type formulae below (cf. Theorem 3.7 and Theorem 3.8). To ensure the finiteness of the Feynman-Kac semigroup, we therefore always assume that $\underline{\mathscr{R}}^E \in \mathrm{K}(M)$ is in the Kato class:

**Definition 3.3** (Kato class). Let $w : M \to \mathbb{C}$ be a measurable function. Then $w$ is said to be in the **contractive Dynkin class** $\mathrm{D}(M)$ (also **extended Kato class**), if there is a $t > 0$ with

$$\sup_{x\in M} \int_0^t \mathbb{E}\left(\mathbb{1}_{\{s<\zeta(x)\}} \left|w(X_s(x))\right|\right) \mathrm{d}s < 1,$$

and $w$ is in the **Kato class** $\mathrm{K}(M)$, if

$$\lim_{t\searrow 0} \sup_{x\in M} \int_0^t \mathbb{E}\left(\mathbb{1}_{\{s<\zeta(x)\}} \left|w(X_s(x))\right|\right) \mathrm{d}s = 0.$$

**Remark 3.4.** (i) The Kato class $\mathrm{K}(M)$ plays an important rôle in the study of Schrödinger operators and their associated semigroups, cf. [Sim82]. The contractive Dynkin class $\mathrm{D}(M)$ appears in [Voi86] to study properties of semigroups associated to Schrödinger operators. In the case of a non-compact manifold it is well-known that there are many technical difficulties with the behaviour of the potentials at $\infty$. The Kato class defines a sufficiently rich class of potentials for which we can still expect the Feynman-Kac formula to make sense pointwise not only vol-a.e. $x$, but *for all* $x \in M$.

(ii) In particular (cf. [Gün17, Remark VI.2.]), in the Euclidean space $\mathbb{R}^m$, we get $\mathrm{L}^q(\mathbb{R}^m) \subset \mathrm{K}(\mathbb{R}^m)$, for $m \geqslant 2$ and $q > \frac{m}{2}$. Then it is well-known, that the Coulomb potential $\frac{1}{|x|}$ is in $\mathrm{K}(\mathbb{R}^3)$.

(iii) Clearly, the classes depend on the Riemannian structure of $M$ and we have

$$\mathrm{K}(M) \subset \mathrm{D}(M).$$

In view of those implications and since it is more common to work with Kato classes, we note that in what follows all assumptions may be relaxed from $\mathrm{K}(M)$ to $\mathrm{D}(M)$.

**Lemma 3.5** ([Gün17, Lemma VI.8.]). *For any $w \in \mathsf{K}(M)$ and $\gamma > 1$, there is $c_\gamma = c_\gamma(w) > 0$, such that*

$$\sup_{x\in M} \mathbb{E}\left(\mathbb{1}_{\{t<\zeta(x)\}} e^{\int_0^t |w(X_s(x))| ds}\right) \leqslant \gamma e^{tc_\gamma} < \infty, \qquad \text{for all } \ t \geqslant 0. \tag{3.5}$$

**Remark 3.6.** The previous Lemma 3.5 can be elaborated in the case of potentials in the Dynkin class (cf. [Gün17, Lemma VI.8.]), namely: For any $w \in \mathsf{D}(M)$ there are $c_k = c_k(w) > 0$, for $k \in \{1, 2\}$, such that

$$\sup_{x\in M} \mathbb{E}\left(\mathbb{1}_{\{t<\zeta(x)\}} e^{\int_0^t |w(X_s(x))| ds}\right) \leqslant c_1 e^{tc_2} < \infty, \qquad \text{for all } \ t \geqslant 0.$$

Following the ideas of [DT01], we will write out the Bismut-type formulae derived for the special cases in Example 3.1 and Example 3.2.

| *From now on, let $E := \bigwedge \mathsf{T}^*M$ and, to shorten notation, we set $\underline{\mathscr{R}} := \underline{\mathscr{R}}^{\bigwedge \mathsf{T}^*M}$.* |
|---|

By [DT01, Theorem 6.1]], we get immediately the following Bismut-type derivative formulae.

**Theorem 3.7.** *Let $\alpha \in \Omega_{\mathsf{L}^2}(M)$ and $\underline{\mathscr{R}}^- \in \mathsf{K}(M)$. Then for any $v \in \bigwedge \mathsf{T}_x M$, we have the following Bismut-type formulae:*

$$\left((\mathbf{d}P_s\alpha)_x, v\right) = -\mathbb{E}\left( /\!\!/_s^{-1}\alpha(X_s(x))\mathbb{1}_{\{s<\zeta(x)\}}, \mathscr{Q}_s \int_0^s \mathscr{Q}_r^{-1}(dB_r \lrcorner\, \mathscr{Q}_r\dot{\ell}_r)\right), \tag{3.6}$$

$$\left((\boldsymbol{\delta}P_s\alpha)_x, v\right) = -\mathbb{E}\left( /\!\!/_s^{-1}\alpha(X_s(x))\mathbb{1}_{\{s<\zeta(x)\}}, \mathscr{Q}_s \int_0^s \mathscr{Q}_r^{-1}(dB_r \wedge \mathscr{Q}_r\dot{\ell}_r)\right), \tag{3.7}$$

*where*

- *$\tau_D < \zeta(x)$ is the first exit time of $X$ from an open relatively compact neighbourhood $D$ of the point $x$,*
- *$dB := /\!\!/^{-1} \circ dX(x)$ is a Brownian motion in $\mathsf{T}_xM$, i.e. the associated anti-development of the Brownian motion $X(x)$,*
- *$(\ell_r)_{r\in[0,s]}$ is any adapted process in $\bigwedge \mathsf{T}_xM$ with absolutely continuous paths such that for for all sufficiently small $\varepsilon > 0$*

$$\mathbb{E}\left(\int_0^{(s-\varepsilon)\wedge\tau_D} \left|\dot{\ell}_r\right|^2 dr\right)^{1/2} < \infty \qquad \text{and} \qquad \ell_0 = v, \quad \ell_r = 0 \qquad \text{for all } \ r \geqslant (s-\varepsilon)\wedge\tau(x,r).$$

*If, in addition, $\alpha \in \Omega_{\mathsf{L}^2}(M)$ is bounded on this neighbourhood, we may take $\varepsilon = 0$.*

Note the remarkable facts in the formula above that curvature only enters locally around a point $x$ and there is no derivative on the righthand side. Moreover, the process $\ell$ may be chosen nicely, as long as it starts in $v$ and is zero when the process $X$ exits from some relatively compact neighbourhood $D$ of $x$.

Next, we prove a covariant Bismut-type formula in the setting of Example 3.2. In the proof, we briefly outline the strategy and refer the reader to Driver and Thalmaier [DT01] studying the theory in full generality, namely on vector bundles, for all the details. Their idea is to define a suitable martingale, say $N_s$, and stay on the *local* martingale level as long as possible. One then shows that $N_s$ is indeed a martingale and takes expectations at times $s = 0$ and $s = t \wedge \tau$, recalling that a *true* martingale has constant expectations thus $\mathbb{E}N_0 = \mathbb{E}N_{s\wedge\tau}$. Note that this method solely involves the geometry and applies *especially* in the case of non-compact manifolds.

**Theorem 3.8.** *Suppose* $\alpha \in \Omega_{L^2}(M)$, $\underline{\mathscr{R}}^- \in \mathsf{K}(M)$ *and* $\xi \in \widetilde{E}_x^* = \mathsf{T}_x M \otimes \bigwedge \mathsf{T}_x M$. *Let* $(\ell_r)_{r\in[0,s]}$ *a bounded adapted process with absolutely continuous paths in* $\bigwedge \mathsf{T}_x M$ *such that* $\mathbb{E}\left(\int_0^{(s-\varepsilon)\wedge\tau_D} |\dot{\ell}_r|^2 \,\mathrm{d}r\right)^{1/2} < \infty$ *and* $\ell_0 = \xi$, $\ell_r = 0$ *for all* $r \geqslant (s-\varepsilon)\wedge\tau_D$ *and for all sufficiently small* $\varepsilon > 0$. *Then,*

$$\left(\nabla \mathrm{e}^{s/2\Delta}\alpha(x), \xi\right) = -\mathbb{E}\left(\mathcal{Q}_s /\!\!/_s^{-1}\alpha(X_s(x))\mathbb{1}_{\{s<\zeta(x)\}}, U_{s\wedge\tau}^{\ell}\right), \tag{3.8}$$

*where*

$$U_s^{\ell} := \int_0^s \mathcal{Q}_r^{-1}(\mathrm{d}B_r \lrcorner \widetilde{\mathcal{Q}}_r \dot{\ell}_r) + \frac{1}{2}\int_0^s (\mathcal{Q}_r)^{-1}\varrho_{/\!/_s}\widetilde{\mathcal{Q}}_r \ell_r \,\mathrm{d}r,$$

$\varrho_{/\!/_s} := /\!\!/_s^{-1}\varrho\, /\!\!/_s$ *with* $\varrho$ *given by* (3.3), $\mathcal{Q}$ *and* $\widetilde{\mathcal{Q}}$ *are defined by* (3.4) *and* $\widetilde{\mathscr{R}}$ *by* (3.2).

**Proof.** According to [DT01, Proposition 3.2 & Theorem 3.7], the process

$$N_r^{\ell} := \left(\widetilde{\mathcal{Q}}_r /\!\!/_r^{-1}\nabla e^{\Delta(s-r)/2}\alpha(X_r(x)), \ell_r\right) - \left(\mathcal{Q}_r /\!\!/_r^{-1} e^{\Delta(s-r)/2}\alpha(X_r(x)), U_r^{\ell}\right)$$

is a local martingale. The existence of the scalar semigroup is provided by semigroup domination and the assumption $\underline{\mathscr{R}}^- \in \mathsf{K}(M)$:

$$\left|P_s^g \alpha(x)\right| \leqslant P_s^{\underline{\mathscr{R}}^-}|\alpha|(x) < \infty.$$

Thus, by the Burkholder-Davis-Gundy inequality, $N_s^{\ell}$ is a true martingale and evaluating at the times $r = 0$ and $r = t\wedge\tau$ and taking expectations, we get by the martingale property $\mathbb{E}N_0^{\ell} = \mathbb{E}N_{s\wedge\tau}^{\ell}$,

$$\left(\nabla \mathrm{e}^{s/2\Delta}\alpha(x), \xi\right) = -\mathbb{E}\left(\mathcal{Q}_{s\wedge\tau} /\!\!/_{s\wedge\tau}^{-1} e^{\Delta(s-s\wedge\tau)/2}\alpha(X_{s\wedge\tau}(x)), U_{s\wedge\tau}^{\ell}\right). \tag{3.9}$$

Note that, by the strong Markov property,

$$\mathcal{Q}_{s\wedge\tau} /\!\!/_{s\wedge\tau}^{-1}\left(P_{s-s\wedge\tau}\alpha\right)(X_{s\wedge\tau}(x)) = \mathbb{E}^{\mathscr{F}_{s\wedge\tau}}\left(\mathcal{Q}_s /\!\!/_s^{-1}\alpha(X_s(x))\mathbb{1}_{\{s<\zeta(x)\}}\right),$$

which is by definition a bounded $\mathscr{F}_{s\wedge\tau}$-measurable random variable. Hence,

$$\left(\nabla \mathrm{e}^{s/2\Delta}\alpha(x), \xi\right) = -\mathbb{E}\left(\mathcal{Q}_s /\!\!/_s^{-1}\alpha(X_s(x))\mathbb{1}_{\{s<\zeta(x)\}}, U_{s\wedge\tau}^{\ell}\right).$$

See also [DT01, p. 82f] and [Tha97] for the technical details. ■

From now on, set $D := \mathsf{B}(x, R)$ to be a ball with small radius, say $R = 1$, and we define

$$\overline{K}(x) := \max\left\{(\mathscr{R}(v), v) : v \in \bigwedge \mathsf{T}_y M, |v| = 1, y \in \mathsf{B}(x,1)\right\}, \tag{3.10}$$

$$\underline{K}(x) := \min\left\{(\mathscr{R}(v), v) : v \in \bigwedge \mathsf{T}_y M, |v| = 1, y \in \mathsf{B}(x,1)\right\}. \tag{3.11}$$

**Theorem 3.9.** *Let $\alpha \in \Omega_{\mathrm{L}^2}(M)$ and $\underline{\mathscr{R}}^- \in \mathsf{K}(M)$. Then, for all $s > 0$,*

$$\left|(\mathbf{d}P_s\alpha)_x\right|^2 \leqslant \Psi(x,s)\Phi(x,s)\,\|\alpha\|^2_{\Omega_{\mathrm{L}^2}(M)}, \tag{3.12}$$

$$\left|(\boldsymbol{\delta}P_s\alpha)_x\right|^2 \leqslant \Psi(x,s)\Phi(x,s)\,\|\alpha\|^2_{\Omega_{\mathrm{L}^2}(M)}, \tag{3.13}$$

*where*

$$\Psi(x,s) := \frac{1}{\sqrt{s}}\exp\left[D(\gamma, c_\gamma(\underline{\mathscr{R}}^-), c_q^{1/q})s + \left(\pi\sqrt{(m-1)\underline{K}(x)^-} + \pi^2(m+5) + (\overline{K}(x)+\underline{K}(x))^-\right)\frac{s}{2}\right], \tag{3.14}$$

*and the finite constant $D$ depends on the constant $c_\gamma(\underline{\mathscr{R}}^-)$ in (3.5) and the constant $c_q$ from the Burkholder-Davis-Gundy inequality.*

**Proof.** We first assume the form $\alpha$ to be real-valued. By Gronwall's inequality, we have

$$|\mathscr{Q}_s|_{\mathrm{op}} \leqslant \exp\left(-\frac{1}{2}\int_0^s \underline{\mathscr{R}}^-(X_r(x))\mathrm{d}r\right) \qquad \text{for all } s \geqslant 0.$$

and hence

$$|\mathscr{Q}_s|_{\mathrm{op}} \leqslant \mathrm{e}^{\underline{K}(x)s/2}, \quad |\mathscr{Q}_s^{-1}|_{\mathrm{op}} \leqslant \mathrm{e}^{\overline{K}(x)s/2} \qquad \mathbb{P}\text{-a.s. on } \{s \leqslant \tau(x,1)\}. \tag{3.15}$$

Let $q \in [2,\infty)$. By a proper choice of the Cameron-Martin space valued process $\ell_s$, it is well-known how to estimate the second factor (cf. [TW98, Proof of Corollary 5.1], [TW11, Remark 3.2], [CTT18]), for some $v \in \bigwedge \mathsf{T}_x M$,

$$|\ell| \leqslant |v|, \qquad \left[\mathbb{E}\left(\int_0^{s\wedge\tau(x,1)} |\dot{\ell}_r|^2\,\mathrm{d}r\right)^q\right]^{1/(2q)} \leqslant \frac{1}{\sqrt{s}}\mathrm{e}^{C(m,2q,R,\underline{K}^-)s/2}\,|v|, \tag{3.16}$$

where

$$C(m,q,R,\underline{K}^-) := \frac{\pi}{2R}\sqrt{(m-1)\underline{K}^-} + \frac{\pi^2}{4R^2}(m+q+3).$$

By the Burkholder-Davis-Gundy inequality, we get

$$\mathbb{E}\left|\int_0^s \mathscr{Q}_r^{-1}(\mathrm{d}B_r \lrcorner\, \mathscr{Q}_r\dot{\ell}_r)\right|^{2q} \leqslant c_q \mathrm{e}^{q(\overline{K}+\underline{K})^- s/2}\,\mathbb{E}\left(\int_0^s |\dot{\ell}_r|^2\,\mathrm{d}r\right)^q. \tag{3.17}$$

By Lemma 3.5, for any $\gamma > 1$, there is a constant $c_\gamma = c_\gamma(\underline{\mathscr{R}}^-)$ such that

$$\sup_{x\in M}\mathbb{E}\left(\mathbb{1}_{\{s<\zeta(x)\}}\mathrm{e}^{\int_0^s|\underline{\mathscr{R}}^-(X_s(x))|\mathrm{d}s}\right) \leqslant \gamma\mathrm{e}^{sc_\gamma} < \infty. \tag{3.18}$$

Now, we can estimate as follows. Let $|v| \leqslant 1$, then using Hölder and Cauchy-Schwarz inequality,

$$\begin{aligned}
\left|(\mathbf{d}P_s\alpha)_x\right| &\leqslant \left[\mathbb{E}\left|\alpha(X_s(x))\mathbb{1}_{\{s<\zeta(x)\}}\right|^p\right]^{1/p}\left[\mathbb{E}\left|\mathscr{Q}_s\mathbb{1}_{\{s<\zeta(x)\}}\int_0^s \mathscr{Q}_r^{-1}(\mathrm{d}B_r \lrcorner\, \mathscr{Q}_r\dot{\ell}_r)\right|^q\right]^{1/q}\\
&\leqslant \left[\mathbb{E}\left|\alpha(X_s(x))\mathbb{1}_{\{s<\zeta(x)\}}\right|^p\right]^{1/p}\left[\mathbb{E}\left(|\mathscr{Q}_s|^{2q}\mathbb{1}_{\{s<\zeta(x)\}}\right)\right]^{1/(2q)}\left[\mathbb{E}\left(\int_0^s \mathscr{Q}_r^{-1}(\mathrm{d}B_r \lrcorner\, \mathscr{Q}_r\dot{\ell}_r)\right)^{2q}\right]^{1/(2q)}\\
&\underset{(3.17)}{\overset{(3.18)}{\leqslant}} \left[\mathbb{E}\left(|\alpha|^p(X_s(x))\mathbb{1}_{\{s<\zeta(x)\}}\right)\right]^{1/p}\gamma\mathrm{e}^{sc_\gamma}c_q^{1/q}\mathrm{e}^{(\overline{K}+\underline{K})^- s/2}\left[\mathbb{E}\left(\int_0^s|\dot{\ell}_r|^2\,\mathrm{d}r\right)^q\right]^{1/(2q)}\\
&\overset{(3.16)}{\leqslant} \mathrm{e}^{C_1(\gamma,c_\gamma,c_q^{1/q})s\,+\,(\overline{K}+\underline{K})^- s/2}\mathrm{e}^{C(m,2q,R,\underline{K}^-)s/2}s^{-1/2}\left[\int_M p_s^{g,(0)}(x,y)\,|\alpha(y)|^p\,\mathrm{vol}_g(\mathrm{d}y)\right]^{1/p}\\
&\leqslant \sqrt{\Psi(x,s)}\sqrt[p]{\Phi(x,s)}\,\|\alpha\|_{\Omega_{\mathrm{L}^p}(M)},
\end{aligned}$$

where $\Psi(x,s)$ is given by (3.14). In particular, for $p=2=q$ the result follows.

By an analogous calculation, we obtain the estimate (3.13). Since complexifications are norm preserving, the proof is complete. ∎

Using similar techniques as in the Proof of the previous theorem, we can show the following estimate. We recall that $\mathsf{R}$ denotes the Riemannian curvature tensor, cf. Example 3.2 above.

**Theorem 3.10.** *Let $\alpha \in \Omega_{\mathsf{L}^2}(M)$ and $\underline{\mathscr{R}}^- \in \mathsf{K}(M)$. Then, for all $\xi \in \mathsf{T}_xM \otimes \bigwedge \mathsf{T}_xM$ and $s>0$,*

$$\left|\left(\nabla P_s\alpha,\xi\right)\right|^2 \leqslant |\xi|^2\,\Xi(x,s)\Phi(x,s)\,\|\alpha\|^2_{\Omega_{\mathsf{L}^2}(M)}, \tag{3.19}$$

*where*

$$\begin{aligned}\Xi(x,s) &:= \mathrm{e}^{D(\gamma,c_\gamma(\underline{\mathscr{R}}^-),c_q^{1/q})s\,+\,(\overline{K}(x)+\underline{K}(x))^-s/2}\left[\mathrm{e}^{C(m,\underline{K}(x)^-)s/2}s^{-1/2}+\max_{y\in\mathsf{B}(x,1)}|\nabla\mathsf{R}(y)|\,s\right]\\ &= \Psi(x,s)+s^{-3/2}\Psi(x,s)\max_{y\in\mathsf{B}(x,1)}|\nabla\mathsf{R}(y)|\,,\end{aligned} \tag{3.20}$$

*and the finite constant $D$ depends on the constant $c_\gamma(\underline{\mathscr{R}}^-)$ in (3.5) and the constant $c_q$ from the Burkholder-Davis-Gundy inequality.*

**Proof.** We first assume the form $\alpha$ to be real-valued. Set

$$U_s^\ell = \int_0^s \mathcal{Q}_r^{-1}(\mathrm{d}B_r \lrcorner\, \widetilde{\mathcal{Q}}_r\dot{\ell}_r) \;+\; \frac{1}{2}\int_0^s \mathcal{Q}_r^{-1}\varrho_{/\!/_r}\widetilde{\mathcal{Q}}_r\ell_r\,\mathrm{d}s = \ell_s^{(1)}+\frac{1}{2}\ell_s^{(2)},$$

where we define processes

$$\ell_s^{(1)} := \int_0^s(\mathcal{Q}_r)^{-1}(\mathrm{d}B_r \lrcorner\, \widetilde{\mathcal{Q}}_r\dot{\ell}_r) \qquad \text{and} \qquad \ell_s^{(2)} := \int_0^s(\mathcal{Q}_r)^{-1}\varrho_{/\!/_r}\widetilde{\mathcal{Q}}_r\ell_r\,\mathrm{d}s.$$

Then $\ell^{(1)}$ is a continuous local martingale and $\ell^{(2)}$ is a continuous process of finite variation. As in the previous Proof of Theorem 3.9, we find

$$|\ell|\leqslant|\xi|\,,\qquad \left[\mathbb{E}\left(\int_0^{s\wedge\tau(x,1)}\left|\dot{\ell}_r\right|^2\mathrm{d}r\right)^q\right]^{1/(2q)} \leqslant \frac{1}{\sqrt{s}}\mathrm{e}^{C(m,2q,R,\underline{K}^-)s/2}\,|\xi|\,.$$

Again using Gronwall's inequality we get (3.15), and by

$$\left|\varrho(X_s(x))\right| \leqslant \max_{y\in\mathsf{B}(x,1)}|\varrho(y)| \leqslant \max_{y\in\mathsf{B}(x,1)}|\nabla\mathsf{R}(y)| \qquad \mathbb{P}-\text{a.s. on }\{s\leqslant\tau(x,1)\}\,, \tag{3.21}$$

we have

$$\mathbb{E}\left|\int_0^s \mathcal{Q}_s^{-1}\varrho_{/\!/_s}\widetilde{\mathcal{Q}}_s\ell_s\,\mathrm{d}s\right| \leqslant e^{(\overline{K}_1+\underline{K}_2)^-s/2}\max_{y\in\mathsf{B}(x,1)}|\nabla\mathsf{R}(y)|\,s\,|\xi|\,. \tag{3.22}$$

By Lemma 3.5, for any $\gamma>1$, there is a constant $c_\gamma=c_\gamma(\underline{\mathscr{R}}^-)$ such that

$$\sup_{x\in M}\mathbb{E}\left(\mathbb{1}_{\{t<\zeta(x)\}}\mathrm{e}^{\int_0^s|\underline{\mathscr{R}}^-(X_s(x))|\mathrm{d}s}\right)\leqslant\gamma\mathrm{e}^{sc_\gamma}<\infty.$$

As in the proof of Theorem 3.9, a similar calculation shows, using Hölder's inequality and the elementary inequality $(a+b)^c\leqslant 2^{c-1}(a^c+b^c)$,

$$\left|\left(\nabla P_s\alpha(x),\xi\right)\right| = \left|\mathbb{E}\left(/\!/_s^{-1}\alpha(X_s(x))\mathbb{1}_{\{s<\zeta(x)\}},\mathcal{Q}_sU^\ell_{s\wedge\tau}\right)\right|$$

$$
\begin{aligned}
&\leqslant 2\left|\mathbb{E}\left(/\!/_s^{-1}\alpha(X_s(x))\mathbb{1}_{\{s<\zeta(x)\}}, \mathcal{Q}_s\ell^{(1)}_{s\wedge\tau}\right)\right| + \left|\mathbb{E}\left(/\!/_s^{-1}\alpha(X_s(x))\mathbb{1}_{\{s<\zeta(x)\}}, \mathcal{Q}_s\ell^{(2)}_{s\wedge\tau}\right)\right| \\
&\leqslant \left[\mathbb{E}\left|\alpha(X_s(x))\mathbb{1}_{\{s<\zeta(x)\}}\right|^p\right]^{1/p}\left(2\left[\mathbb{E}\left(\mathcal{Q}_s\mathbb{1}_{\{s<\zeta(x)\}}\ell^{(1)}_{s\wedge\tau}\right)^q\right]^{1/q} + \left[\mathbb{E}\left(\mathcal{Q}_s\mathbb{1}_{\{s<\zeta(x)\}}\ell^{(2)}_{s\wedge\tau}\right)^q\right]^{1/q}\right) \\
&\leqslant |\xi|\left[\mathbb{E}\left(|\alpha|^p(X_s(x))\mathbb{1}_{\{s<\zeta(x)\}}\right)\right]^{1/p}\gamma e^{sc_\gamma}c_q^{1/q}e^{(\overline{K}+\underline{K})^- s/2}\left(2\left[\mathbb{E}\left(\int_0^s\left|\dot\ell_r\right|^2 dr\right)^q\right]^{1/(2q)} + \right. \\
&\qquad\qquad \left. + s\max_{y\in B(x,1)}|\nabla R(y)|\right) \\
&\leqslant |\xi|\, e^{D(\gamma,c_\gamma,c_q^{1/q})s\,+\,(\overline{K}+\underline{K})^- s/2}\left[e^{C(m,2q,R,\underline{K}^-)s/2}s^{-1/2} + \left(\max_{y\in B(x,1)}|\nabla R(y)|\right)s\right]\times \\
&\qquad\qquad \times\left[\int_M p_s^{g,(0)}(x,y)\,|\alpha(y)|^p\,\mathrm{vol}_g(dy)\right]^{1/p} \\
&\leqslant |\xi|\sqrt{\Xi(x,s)}\,\sqrt[p]{\Phi(x,s)}\,\|\alpha\|_{\Omega_{L^p}(M)}\,.
\end{aligned}
$$

where $\Xi(x,s)$ is given by (3.20).

In particular, for $p = 2 = q$ the result follows. Since complexification is norm preserving, the Proof is complete. ■

In the remaining part of this section, we will make use of Theorem 3.10 and derive similar estimates for the exterior derivative and codifferential transformed by the smooth vector bundle morphism $\mathcal{A}_{g,h}$ defined in (2.4). The key observation will be Proposition 3.15 showing how to estimate the transformed codifferential $\boldsymbol{\delta}_g$ (with respect the metric $g$) applied to the semigroup in terms of covariant derivative $\nabla^g$ of $g$ applied to the semigroup. In addition, using Lemma 2.8, a direct consequence is an analogous result (cf. Corollary 3.18) in terms of the new metric, i.e. for the transformed codifferential $\boldsymbol{\delta}_h$ (with respect the metric $h$) applied to the semigroup.

To this end, we will make use of the well-known metric descriptions of the exterior derivative $\mathbf{d}$ and the codifferential $\boldsymbol{\delta}_g$ (with respect the metric $g$), namely:

**Lemma 3.11** ([Jos17, Lemma 4.3.4]). *Let $e_1, \dots, e_m$ (where $m = \dim M$) be a local orthonormal frame for $T_xM$ ($x \in M$) and $\varepsilon^1, \dots, \varepsilon^m$ be its dual coframe, i.e. $\varepsilon^j(e_i) = \delta_i^j$. Then*

$$
\mathbf{d} = \sum_{i=1}^m \varepsilon^i \wedge \nabla^g_{e_i} \quad \text{and} \quad \boldsymbol{\delta}_g = -\sum_{i=1}^m e_i \lrcorner \nabla^g_{e_i}. \tag{3.23}
$$

For the remainder of this section, let $g, h \in \mathrm{Metr}M$, $g \sim h$, and $\mathcal{A} := \mathcal{A}_{g,h}$ be the smooth vector bundle morphism defined by (2.4).

**Proposition 3.12.** *Let $\alpha \in \Omega_{L^2}(M)$ and $\underline{\mathcal{R}}^- \in K(M)$. Then, for any orthonormal frame $(e_i)_{i=1}^m$ for $T_xM$ ($x \in M$) and dual coframe $(\varepsilon^i)_{i=1}^m$, we decompose*

$$
\mathbf{d}(\mathcal{A}^{1/2}P_s\alpha) = \sum_{i=1}^m \left(\varepsilon^i \wedge \nabla^g_{e_i}P_s\alpha \circ \mathcal{A}^{1/2} + \varepsilon^i \wedge P_s\alpha \circ \left(\nabla^g_{e_i}\mathcal{A}^{1/2}\right)\right).
$$

**Proof.** From (3.23), we get for any $\eta \in \Omega(M)$

$$\mathbf{d}(\mathscr{A}^{1/2}\eta) = \mathbf{d}(\eta \circ \mathscr{A}^{1/2}) = \sum_{i=1}^{m} \varepsilon^i \wedge \nabla^g_{e_i}\eta \circ \mathscr{A}^{1/2} + \sum_{i=1}^{m} \varepsilon^i \wedge \eta \circ \nabla^g_{e_i}\mathscr{A}^{1/2},$$

where $(e_i)_{i=1}^m$ denotes an orthonormal frame for $\mathrm{T}_xM$ $(x \in M)$ and $(\varepsilon^i)_{i=1}^m$ the dual coframe. In particular, for $\eta = P_s\alpha$ the claim follows. ∎

**Corollary 3.13.** *Let $\alpha \in \Omega_{\mathrm{L}^2}(M)$ and $\underline{\mathscr{R}}^- \in \mathrm{K}(M)$. Then,*

$$\left|\mathscr{A}^{-1/2}\mathbf{d}(\mathscr{A}^{1/2}P_s\alpha(x))\right|^2 \lesssim \left(\delta^{\nabla}_{g,h}(x) + \Xi(x,s)\right)\Phi(x,s)\,\|\alpha\|^2_{\Omega_{\mathrm{L}^2}(M)}\,. \tag{3.24}$$

**Proof.** By Theorem 3.10, we have

$$\left|\nabla^g_{e_i}P_s\alpha \circ \mathscr{A}^{1/2}\right|^2 \leqslant C(m)\left|\mathscr{A}^{1/2}(x)\right|^2 \Xi(x,s)\Phi(x,s)\,\|\alpha\|^2_{\Omega_{\mathrm{L}^2}(M)}$$

and, combined with Lemma 3.5 and Lemma 2.4,

$$\begin{aligned}\left|P_s\alpha \circ \left(\nabla^g_{e_i}\mathscr{A}^{1/2}\right)\right|^2 &\leqslant C(m,\gamma,c_\gamma,s)\left|\nabla^g\mathscr{A}^{1/2}(x)\right|^2_g \Phi(x,s)\,\|\alpha\|^2_{\Omega_{\mathrm{L}^2}(M)}\\ &\lesssim \left|\mathscr{A}^{1/2}(x)\right|^2 \delta^{\nabla}_{g,h}(x)\Phi(x,s)\,\|\alpha\|^2_{\Omega_{\mathrm{L}^2}(M)}\,.\end{aligned}$$

Thus the claim follows:

$$\left|\mathscr{A}^{-1/2}\mathbf{d}(\mathscr{A}^{1/2}P_s\alpha(x))\right|^2 \lesssim \left|\mathscr{A}^{-1/2}(x)\right|^2\left|\mathscr{A}^{1/2}(x)\right|^2\left(\delta^{\nabla}_{g,h}(x) + \Xi(x,s)\right)\Phi(x,s)\,\|\alpha\|^2_{\Omega_{\mathrm{L}^2}(M)}\,.$$ ∎

**Remark 3.14.** By a similar calculation the estimate also holds if we interchange the rôles of $\mathscr{A}^{1/2}$ and $\mathscr{A}^{-1/2}$ in the previous Corollary 3.13.

**Proposition 3.15.** *Let $\alpha \in \Omega_{\mathrm{L}^2}(M)$ and $\underline{\mathscr{R}}^- \in \mathrm{K}(M)$. Then, for any orthonormal basis $(e_i)_i$ for $\mathrm{T}_xM$ $(x \in M)$, we decompose*

$$\boldsymbol{\delta}_g(\mathscr{A}^{1/2}P_s\alpha) = -\sum_{i=1}^{m} e_i \lrcorner \left(\nabla^g_{e_i}P_s\alpha\right)\circ \mathscr{A}^{1/2} - \sum_{i=1}^{m} e_i \lrcorner P_s\alpha \circ \left(\nabla^g_{e_i}\mathscr{A}^{1/2}\right). \tag{3.25}$$

**Proof.** Let $\eta \in \Omega(M)$ be arbitrary and $(e_i)_{i=1}^m$ an orthonormal basis for $\mathrm{T}_xM$ $(x \in M)$. By definition of the codifferential, (3.23) in Lemma 3.11 above, we get

$$\begin{aligned}\boldsymbol{\delta}_g(\mathscr{A}^{1/2}\eta) &= -\sum_{i=1}^{m} e_i \lrcorner \nabla^g_{e_i}(\mathscr{A}^{1/2}\eta)\\ &= -\sum_{i=1}^{m}\left(e_i \lrcorner (\mathscr{A}^{1/2}\nabla^g_{e_i}\eta) + e_i \lrcorner ((\nabla^g_{e_i}\mathscr{A}^{1/2})\eta)\right).\end{aligned}$$

In particular, if we set $\eta = P_s\alpha(x)$, then the following equalities hold:

$$\sum_{i=1}^{m} e_i \lrcorner \left(\mathscr{A}^{1/2}\nabla^g_{e_i}P_s\alpha\right) = \sum_{i=1}^{m} e_i \lrcorner \left(\nabla^g_{e_i}P_s\alpha\right)\circ \mathscr{A}^{1/2}$$

and

$$\sum_{i=1}^{m} e_i \lrcorner ((\nabla^g_{e_i}\mathscr{A}^{1/2})P_s\alpha) = \sum_{i=1}^{m} e_i \lrcorner P_s\alpha \circ \left(\nabla^g_{e_i}\mathscr{A}^{1/2}\right).$$ ∎

**Corollary 3.16.** *Let $\alpha \in \Omega_{\mathrm{L}^2}(M)$ and $\underline{\mathscr{R}}^- \in \mathrm{K}(M)$. Then,*

$$\left|\mathscr{A}^{-1/2}\boldsymbol{\delta}_g(\mathscr{A}^{1/2}P_s\alpha(x))\right|^2 \lesssim \left(\delta^{\nabla}_{g,h}(x) + \Xi(x,s)\right)\Phi(x,s)\,\|\alpha\|^2_{\Omega_{\mathrm{L}^2}(M)}. \tag{3.26}$$

**Proof.** By Theorem 3.10 we have

$$\left|\sum_{i=1}^m e_i \lrcorner \left(\nabla^g_{e_i}P_s\alpha\right)\circ\mathscr{A}^{1/2}\right|^2 \leqslant C(m)\left|\mathscr{A}^{1/2}(x)\right|^2 \Xi(x,s)\Phi(x,s)\,\|\alpha\|^2_{\Omega_{\mathrm{L}^2}(M)}$$

and, combined with Lemma 3.5 and Lemma 2.4,

$$\begin{aligned}\left|e_i \lrcorner P_s\alpha\circ\left(\nabla^g_{e_i}\mathscr{A}^{1/2}\right)\right|^2 &\leqslant C(m,\gamma,c_\gamma,s)\left|\nabla^g\mathscr{A}^{1/2}(x)\right|^2_g\Phi(x,s)\,\|\alpha\|^2_{\Omega_{\mathrm{L}^2}(M)}\\ &\lesssim \left|\mathscr{A}^{1/2}(x)\right|^2\delta^{\nabla}_{g,h}(x)\Phi(x,s)\,\|\alpha\|^2_{\Omega_{\mathrm{L}^2}(M)}.\end{aligned}$$

Thus the claim follows:

$$\left|\mathscr{A}^{-1/2}\boldsymbol{\delta}_g(\mathscr{A}^{1/2}P_s\alpha(x))\right|^2 \lesssim \left|\mathscr{A}^{-1/2}(x)\right|^2\left|\mathscr{A}^{1/2}(x)\right|^2\left(\delta^{\nabla}_{g,h}(x) + \Xi(x,s)\right)\Phi(x,s)\,\|\alpha\|^2_{\Omega_{\mathrm{L}^2}(M)}. \qquad \blacksquare$$

**Remark 3.17.** By a similar calculation the estimate also holds if we interchange the rôles of $\mathscr{A}^{1/2}$ and $\mathscr{A}^{-1/2}$ in the previous Corollary 3.16.

Finally, we obtain a similar estimate for the transformed codifferential $\boldsymbol{\delta}_h$ (with respect the metric $h$) applied to the semigroup. As we can express $\boldsymbol{\delta}_h$ in terms of $\boldsymbol{\delta}_g$ using Lemma 2.8, we want to emphasise that the fibrewise norm and the involved quantities are taken with respect to the metric $g$.

**Corollary 3.18.** *Let $\alpha \in \Omega_{\mathrm{L}^2}(M)$ and $\underline{\mathscr{R}}^- \in \mathrm{K}(M)$. Then,*

$$\left|\mathscr{A}^{1/2}\boldsymbol{\delta}_h(\mathscr{A}^{-1/2}P_s\alpha(x))\right|^2_g \lesssim \left(\delta^{\nabla}_{g,h}(x) + \Xi_g(x,s)\right)\Phi_g(x,s)\,\|\alpha\|^2_{\Omega_{\mathrm{L}^2}(M,g)}. \tag{3.27}$$

**Proof.** By Lemma 2.8, recall that

$$\boldsymbol{\delta}_h(\eta) = \mathscr{A}^{-1}\left(\boldsymbol{\delta}_g(\mathscr{A}\eta) - (\mathbf{d}\log\rho)^{\sharp_g}\lrcorner_g(\mathscr{A}\eta)\right) \qquad \text{for all } \eta\in\Omega_{\mathrm{C}^\infty}(M),$$

so that

$$(\mathscr{A}^{1/2}\boldsymbol{\delta}_h\mathscr{A}^{-1/2})(\eta) = \mathscr{A}^{-1/2}\left(\boldsymbol{\delta}_g(\mathscr{A}^{1/2}\alpha) - (\mathbf{d}\log\rho)^{\sharp_g}\lrcorner_g(\mathscr{A}^{1/2}\eta)\right).$$

Using Proposition 2.12 combined with Lemma 2.5, in addition

$$\begin{aligned}\mathscr{A}^{-1/2}\left|(\mathbf{d}\log\rho)^{\sharp_g}\lrcorner(\mathscr{A}^{1/2}\eta)\right|_g &\leqslant C(m)\left|\mathscr{A}^{-1/2}\right|\left|\mathscr{A}^{-1}\right|\left|\nabla^g\mathscr{A}\right|_g\left|\mathscr{A}^{1/2}\right|\left|\eta\right|_g\\ &\leqslant C(m)\left|\mathscr{A}^{-1/2}\right|\left|\mathscr{A}^{-1}\right|\left|\mathscr{A}\right|\left|\nabla^h_X - \nabla^g_X\right|_g\left|\mathscr{A}^{1/2}\right|\left|\eta\right|_g\\ &\leqslant C(m)\left|\nabla^h_X - \nabla^g_X\right|_g\left|\eta\right|_g.\end{aligned}$$

Finally, combined with Corollary 3.16,

$$\left|\mathscr{A}^{1/2}\boldsymbol{\delta}_h(\mathscr{A}^{-1/2}P_s\alpha(x))\right|^2_g \lesssim \left(\delta^{\nabla}_{g,h}(x) + \Xi_g(x,s)\right)\Phi_g(x,s)\,\|\alpha\|^2_{\Omega_{\mathrm{L}^2}(M,g)}$$

which proves the claim. $\blacksquare$

# 4 Main results

Since we assume $M$ to be geodesically complete, we can restrict ourselves to smooth compactly supported differential forms. Using the common abuse of notation, the unique realisations of the exterior derivative $\mathbf{d}$, the codifferential $\boldsymbol{\delta}_g$ and Hodge Laplacian $\Delta_g$ will be denoted by the same symbol.

In addition, we define the operators

$$(\hat{P}_s^g)_{s>0} := (\mathbf{d}\, P_s^g)_{s>0} \subset \mathscr{L}(\Omega_{\mathsf{L}^2}(M,g), \Omega_{\mathsf{L}^2}(M,g)),$$
$$(\check{P}_s^g)_{s>0} := (\boldsymbol{\delta}_g P_s^g)_{s>0} \subset \mathscr{L}(\Omega_{\mathsf{L}^2}(M,g), \Omega_{\mathsf{L}^2}(M,g)),$$

$$(\hat{P}_s^{g,h})_{s>0} := (I_{g,h}^{-1}\mathbf{d}\, I_{g,h} P_s^g)_{s>0} \subset \mathscr{L}(\Omega_{\mathsf{L}^2}(M,g), \Omega_{\mathsf{L}^2}(M,g)),$$
$$(\check{P}_s^{g,h})_{s>0} := (I_{g,h}\boldsymbol{\delta}_g I_{g,h}^{-1} P_s^h)_{s>0} \subset \mathscr{L}(\Omega_{\mathsf{L}^2}(M,g), \Omega_{\mathsf{L}^2}(M,g)).$$

Let $\hat{p}_s^g(x,y)$, $\check{p}_s^g(x,y)$, $\hat{p}_s^{g,h}(x,y)$ and $\check{p}_s^{g,h}(x,y)$ be the corresponding jointly smooth integral kernel of $\hat{P}_s^g$, $\check{P}_s^g$, $\hat{P}_s^{g,h}$ and $\check{P}_s^{g,h}$, respectively. For example, recall that this implies [Gün17, Theorem II.1.] that

$$(0,\infty) \times M \times M \ni (s,x,y) \mapsto \hat{p}_s(x,y) \in \mathsf{Hom}\left(\bigwedge \mathsf{T}_y^* M, \bigwedge \mathsf{T}_x^* M\right)$$

is the uniquely determined map such that we have

$$\hat{P}_s^g \alpha(x) = \int_M \hat{p}_s^g(x,y)\alpha(y)\,\mathrm{vol}_g(\mathrm{d}y) \qquad \text{for all } s>0,\ \alpha \in \Omega_{\mathsf{L}^2}(M,g),\ x \in M.$$

By Riesz' representation theorem, the next result follows from the gradient estimates, Theorem 3.9, for the exterior derivative and the codifferential.

**Theorem 4.1.** *For every $g \in \mathsf{Metr}M$, $(s,x) \in (0,\infty) \times M$, we have*

$$\int \left|\hat{p}_s^g(x,y)\right|_g^2 \mathrm{vol}_g(\mathrm{d}y) \leqslant \Psi_g(x,s)\Phi_g(x,s), \tag{4.1}$$
$$\int \left|\check{p}_s^g(x,y)\right|_g^2 \mathrm{vol}_g(\mathrm{d}y) \leqslant \Psi_g(x,s)\Phi_g(x,s). \tag{4.2}$$

By Riesz' representation theorem, the next result follows from the estimates in Corollaries 3.13, 3.16 and 3.18 for the transformed exterior derivative and for the transformed codifferential with respect to $g$, and $h$, respectively.

**Theorem 4.2.** *For every $g \in \mathsf{Metr}M$, $(s,x) \in (0,\infty) \times M$, we have*

$$\int \left|\hat{p}_s^{g,h}(x,y)\right|_g^2 \mathrm{vol}_g(\mathrm{d}y) \lesssim \left(\delta_{g,h}^{\nabla}(x) + \Xi_g(x,s)\right)\Phi_g(x,s), \tag{4.3}$$
$$\int \left|\check{p}_s^{g,h}(x,y)\right|_g^2 \mathrm{vol}_g(\mathrm{d}y) \lesssim \left(\delta_{g,h}^{\nabla}(x) + \Xi_g(x,s)\right)\Phi_g(x,s), \tag{4.4}$$
$$\int \left|\check{p}_s^{h,g}(x,y)\right|_g^2 \mathrm{vol}_g(\mathrm{d}y) \lesssim \left(\delta_{g,h}^{\nabla}(x) + \Xi_g(x,s)\right)\Phi_g(x,s). \tag{4.5}$$

We can now state the main result on the existence and completeness of the wave operators $W_\pm(\Delta_h, \Delta_g, I)$ implying the corresponding spectra to coincide.

**Theorem 4.3.** *Let $g, h \in \mathrm{Metr}M$, $g \sim h$, and assume that there exists $C < \infty$ such that $\left|\delta^{\nabla}_{g,h}\right| \leqslant C$ and for both $\nu \in \{g, h\}$, we have $\left|\mathscr{R}_\nu\right|_\nu \in \mathsf{K}(M)$ and*

$$\int \max\left\{\delta_{g,h}(x), \delta^{\nabla}_{g,h}(x) + \Xi_g(x,s), \Psi_\nu(x,s)\right\} \Phi_\nu(x,s) \operatorname{vol}_\nu(\mathrm{d}x) < \infty, \qquad \text{some } s > 0. \tag{4.6}$$

*Then the wave operators*

$$W_\pm(\Delta_h, \Delta_g, I_{g,h}) = \operatorname*{s-lim}_{t\to\pm\infty} \mathrm{e}^{it\Delta_h} I_{g,h} \mathrm{e}^{-it\Delta_g} \mathsf{P}_{\mathrm{ac}}(\Delta_g)$$

*exist and are complete. Moreover, $W_\pm(\Delta_h, \Delta_g, I_{g,h})$ are partial isometries with initial space $\operatorname{ran} \mathsf{P}_{\mathrm{ac}}(\Delta_g)$ and final space $\operatorname{ran} \mathsf{P}_{\mathrm{ac}}(\Delta_h)$, and we have $\sigma_{\mathrm{ac}}(\Delta_g) = \sigma_{\mathrm{ac}}(\Delta_g)$.*

The proof of Theorem 4.3 will be given in Section 5.

In the special case of for 0-forms, i.e. functions, and the Hodge Laplacian acting on 0-forms is the Laplace-Beltrami operator. Recall that the Weitzenböck curvature endomorphism $\mathscr{R}^{(k)}$ on 1-forms is given by the Ricci curvature, $\mathscr{R}^{(1),\mathrm{tr}} = \mathrm{Ric}$. Then we get the following result similar to the main result of [GT20, Theorem A].

**Corollary 4.4.** *Let $g, h \in \mathrm{Metr}M$, $g \sim h$, and assume that the function $\delta^{\nabla}_{g,h}$ is bounded, and for some $s \in (0, \infty)$ and both $\nu \in \{g, h\}$ satisfy (4.6). Let $-\Delta_\nu \geqslant 0$ be the unique self-adjoint extensions of the Laplace-Beltrami operator for $\nu \in \{g, h\}$. Then the wave operators $W_\pm(-\Delta_h, -\Delta_g, I)$ exist and are complete. Moreover, $W_\pm(-\Delta_h, -\Delta_g, I)$ are partial isometries with initial space $\operatorname{ran} \mathsf{P}_{\mathrm{ac}}(-\Delta_g)$ and final space $\operatorname{ran} \mathsf{P}_{\mathrm{ac}}(-\Delta_h)$, and we have $\sigma_{\mathrm{ac}}(-\Delta_h) = \sigma_{\mathrm{ac}}(-\Delta_g)$.*

## 5 Proof of the main result

The next Lemma shows assumption (2) in the Belopol'skii-Birman theorem A.3. As we will see in its proof, it therefore necessary for the *potentials* $\mathscr{R}_\nu \in \mathsf{K}(M)$ to be in the Kato class, not only $\underline{\mathscr{R}}^-_\nu \in \mathsf{K}(M)$, for $\nu \in \{g, h\}$.

We denote by $\mathbf{q}_\nu$ the nonnegative closed sesquilinear form corresponding to $\Delta_\nu$, i.e. $\mathbf{q}_\nu(\alpha) = \langle \Delta_\nu \alpha, \alpha \rangle = \left\|\mathsf{D}_\nu \alpha\right\|^2_\nu$ with $\operatorname{dom} \mathbf{q}_\nu = \operatorname{dom} \sqrt{\Delta_\nu}$ for any $\nu \in \{g, h\}$.

**Lemma 5.1.** *Let $g, h \in \mathrm{Metr}M$, $g \sim h$, and assume that there exists $C < \infty$ such that $\left|\delta^{\nabla}_{g,h}\right| \leqslant C$ and for both $\nu \in \{g, h\}$, we have $\left|\mathscr{R}_\nu\right|_\nu \in \mathsf{K}(M)$. Then*

$$I_{g,h} \operatorname{dom} \mathbf{q}_g = \operatorname{dom} \mathbf{q}_h.$$

**Proof.** Note that, for any $\nu \in \{g, h\}$, $\operatorname{dom} \mathbf{q}_\nu$ is the closure of compactly supported forms $\Omega_{\mathsf{C}^\infty_c}(M, \nu)$ with respect to the Dirac graph norm

$$\alpha \mapsto \left(\|\alpha\|^2_\nu + \left\|\mathsf{D}_\nu \alpha\right\|^2_\nu\right)^{1/2}.$$

Moreover let $\mathbf{q}^\nabla$ be the nonnegative closed sesquilinear form corresponding to the horizontal Laplacian $\square_\nu = -(\nabla^\nu)^* \nabla^\nu$, i.e. $\mathbf{q}^\nabla_\nu(\alpha) = \langle \square_\nu \alpha, \alpha \rangle_\nu = \|\nabla^\nu \alpha\|^2_\nu$ with $\operatorname{dom}(\mathbf{q}^\nabla_\nu)$ its natural domain of

definition. Recall that by the Weitzenböck formula, we thus have the relation

$$\Delta_\nu = \Box_\nu - \mathscr{R}_\nu .$$

As the Weitzenböck curvature term $\mathscr{R}_\nu$ is in the Kato class, by [Gün17, Lemma VII.4.] the corresponding form domains

$$\operatorname{dom} \mathbf{q}_\nu = \operatorname{dom} \mathbf{q}_\nu^\nabla$$

coincide. To this end, it suffices to show that

$$I \operatorname{dom} \mathbf{q}_g^\nabla = \operatorname{dom} \mathbf{q}_h^\nabla .$$

For all compactly supported $\alpha \in \Omega_{C_c^\infty}(M,g)$, we write

$$\begin{aligned}\nabla^h(I\alpha) &= \nabla^h(\mathscr{A}^{-1/2}\alpha) = \nabla^h\mathscr{A}^{-1/2}(\alpha) + \mathscr{A}^{-1/2}\nabla^h\alpha \\ &= \left(\nabla^h - \nabla^g\right)\mathscr{A}^{-1/2}(\alpha) + \nabla^g\mathscr{A}^{-1/2}(\alpha) + \mathscr{A}^{-1/2}\left(\nabla^h - \nabla^g\right)\alpha + \mathscr{A}^{-1/2}\nabla^g\alpha .\end{aligned}$$

Moreover,

$$\begin{aligned}\left|\left(\nabla_X^h - \nabla_X^g\right)\mathscr{A}^{-1/2}\right|_h &= \left|\mathscr{A}^{1/2}\left(\nabla_X^h - \nabla_X^g\right)\mathscr{A}^{-1/2}\right|_g \\ &\leqslant \left|\mathscr{A}^{1/2}\right|\left|\nabla_X^h - \nabla_X^g\right|_g\left|\mathscr{A}^{-1/2}\right|\end{aligned}$$

and

$$\left|\mathscr{A}^{-1/2}\left(\nabla_X^h - \nabla_X^g\right)\right|_h = \left|\nabla_X^h - \nabla_X^g\right|_g .$$

Thus we can estimate as follows, with a constant $D < \infty$ that might change from line to line,

$$\begin{aligned}\left\|\nabla^h(I\alpha)\right\|_h^2 &= \int_M \left|\nabla^h(\mathscr{A}^{-1/2}\alpha)\right|_h^2 \mathrm{d}\,\mathrm{vol}_h \\ &= \int_M \left|\left(\nabla^h - \nabla^g\right)\mathscr{A}^{-1/2}(\alpha) + \nabla^g\mathscr{A}^{-1/2}(\alpha) + \mathscr{A}^{-1/2}\left(\nabla^h - \nabla^g\right)(\alpha) + \mathscr{A}^{-1/2}\nabla^g\alpha\right|_h^2 \mathrm{d}\,\mathrm{vol}_h \\ &\leqslant D\int_M \Big(\left|\left(\nabla^h - \nabla^g\right)\mathscr{A}^{-1/2}(\alpha)\right|_h^2 + \left|\nabla^g\mathscr{A}^{-1/2}(\alpha)\right|_h^2 + \\ &\qquad \left|\mathscr{A}^{-1/2}\left(\nabla^h - \nabla^g\right)(\alpha)\right|_h^2 + \left|\mathscr{A}^{-1/2}\nabla^g\alpha\right|_h^2\Big)\rho_{g,h}\, \mathrm{d}\,\mathrm{vol}_g \\ &\leqslant D\int_M \Big(\left|\mathscr{A}^{1/2}\left(\nabla^h - \nabla^g\right)\mathscr{A}^{-1/2}(\alpha)\right|_g^2 + \left|\mathscr{A}^{1/2}\nabla^g\mathscr{A}^{-1/2}(\alpha)\right|_g^2 + \\ &\qquad \left|\left(\nabla^h - \nabla^g\right)(\alpha)\right|_g^2 + \left|\nabla^g\alpha\right|_g^2\Big)\mathrm{d}\,\mathrm{vol}_g \\ &\leqslant D\int_M \Big(\left|\mathscr{A}^{1/2}\right|^2\left|\nabla^h - \nabla^g\right|_g^2\left|\mathscr{A}^{-1/2}\right|^2|\alpha|_g^2 + \left|\mathscr{A}^{1/2}\right|^2|\mathscr{A}|^2\left|\mathscr{A}^{-1/2}\right|^6\left|\nabla_X^h - \nabla_X^g\right|_g^2|\alpha|_g^2 + \\ &\qquad \left|\nabla^h - \nabla^g\right|_g^2|\alpha|_g^2 + |\nabla^g\alpha|_g^2\Big)\mathrm{d}\,\mathrm{vol}_g \\ &\leqslant D\int_M \left(\left|\nabla^h - \nabla^g\right|_g^2|\alpha|_g^2 + |\nabla^g\alpha|_g^2\right)\mathrm{d}\,\mathrm{vol}_g\end{aligned}$$

$$\leqslant D \int_M \left( \left\| \delta^{\nabla}_{g,h} \right\|_\infty |\alpha|_g^2 + |\nabla^g \alpha|_g^2 \right) \mathrm{d}\,\mathrm{vol}_g$$
$$\leqslant D \left( \|\alpha\|_g^2 + \|\nabla^g \alpha\|_g^2 \right),$$

using the elementary inequality $(a+b)^c \leqslant 2^{c-1}(a^c + b^c)$ multiple times and that $\delta^{\nabla}$ is bounded by assumption.

Hence, we arrive at the estimate

$$\|I\alpha\|_h^2 + \left\|\nabla^h I\alpha\right\|_h^2 \leqslant C \left( \|\alpha\|_g^2 + \|\nabla^g \alpha\|_g^2 \right),$$

proving

$$I \operatorname{dom} \mathbf{q}_g \subset \operatorname{dom} \mathbf{q}_h.$$

Since $I^{-1} = I_{g,h}^{-1} = I_{h,g}$ and the arguments above are symmetric in $g$ and $h$, this shows the claim. ■

Next, we denote by $|\cdot| : \mathbb{C} \to \mathbb{R}$ the absolute value function and by $\operatorname{sgn} : \mathbb{C} \to \mathbb{C}$ the sign-function with $\operatorname{sgn}(0) = 1$. We note that [Wei80] if $P$ is normal operator (e.g. positive or diagonalisable), we get the (pointwise) polar decomposition $P = |P| (\operatorname{sgn} P)$, where $|P|(x) = |P(x)| \geqslant 0$ and $|\operatorname{sgn} P(x)| = 1$, and where $|P|(x)$ is a non-negative endomorphism and $\operatorname{sgn} P(x)$ is unitary. More precisely, by the spectral theorem choosing $f(\lambda) := |\lambda|$, we have an endomophism

$$f(P(x)) : \bigwedge \mathsf{T}^* M \to \bigwedge \mathsf{T}^* M$$

giving rise to a decomposition

$$P(x) = |P(x)| \operatorname{sgn} P(x) : \bigwedge \mathsf{T}^* M \to \bigwedge \mathsf{T}^* M.$$

For the proof of Theorem 4.3, we now introduce sections

$$\mathsf{S}_{g,h} : M \to \mathbb{R}$$
$$\mathsf{S}_{g,h}(x) := \rho_{g,h}(x)^{1/2} - \rho_{g,h}(x)^{-1/2} = 2 \sinh \frac{1}{2} \log(\rho_{g,h}(x)),$$

$$\hat{\mathsf{S}}_{g,h} : M \to \operatorname{End}\left(\bigwedge \mathsf{T}^* M\right)$$
$$\hat{\mathsf{S}}_{g,h}(x) := \left(\rho_{g,h}(x)\mathscr{A}_{g,h}(x)\right)^{1/2} - \left(\rho_{g,h}(x)\mathscr{A}_{g,h}(x)\right)^{-1/2} = 2 \sinh \frac{1}{2} \log(\rho_{g,h}(x)\mathscr{A}_{g,h}(x)),$$

$$\hat{\mathsf{S}}_{g,h;\nu} : \Omega_{\mathsf{L}^2}(M, \nu) \to \Omega_{\mathsf{L}^2}(M, \nu)$$
$$\hat{\mathsf{S}}_{g,h;\nu}\alpha(x) := \left|\hat{\mathsf{S}}_{g,h}(x)\right|^{1/2} \alpha(x),$$

$$\mathsf{U}_{g,h} : \Omega_{\mathsf{L}^2}(M, g) \to \Omega_{\mathsf{L}^2}(M, h)$$
$$\mathsf{U}_{g,h}\alpha(x) := \mathscr{A}_{g,h}(x)^{-1/2}\alpha(x),$$

$$\hat{\mathsf{U}}_{g,h} : \Omega_{\mathsf{L}^2}(M, g) \to \Omega_{\mathsf{L}^2}(M, h)$$
$$\hat{\mathsf{U}}_{g,h}\alpha(x) := (\rho_{g,h}(x)\mathscr{A}_{g,h}(x))^{-1/2}\alpha(x),$$

$$\widetilde{\mathsf{U}}_{g,h} : \Omega_{\mathsf{L}^2}(M, g) \to \Omega_{\mathsf{L}^2}(M, h)$$
$$\widetilde{\mathsf{U}}_{g,h}\alpha(x) := (\operatorname{sgn} \hat{\mathsf{S}}_{g,h}(x))(\rho_{g,h}(x)\mathscr{A}_{g,h}(x))^{-1/2}\alpha(x).$$

By quasi-isometry, $g \sim h$, the operators $\hat{\mathbf{S}}_{g,h;\nu}$, $\mathbf{U}_{g,h}$, $\hat{\mathbf{U}}_{g,h}$ and $\widetilde{\mathbf{U}}_{g,h}$ are bounded. Moreover, we get the pointwise estimate [HPW14, Lemma 3.3], [GT20, (4.1)].

**Lemma 5.2.** *We have the pointwise estimate*

$$\max\left\{\left|\mathbf{S}_{g,h}(x)\right|, \sigma_{\max}\left(\left|\hat{\mathbf{S}}_{g,h}(x)\right|\right)\right\} \leqslant \delta_{g,h}(x) \qquad \text{for all } x \in M. \tag{5.1}$$

**Proof.** We write $\rho = \rho_{g,h}$ and $\mathscr{A} = \mathscr{A}_{g,h}$ for short. By definition, we have

$$\left|\hat{\mathbf{S}}_{g,h}\right| = \left|(\rho\mathscr{A})^{1/2} - (\rho\mathscr{A})^{-1/2}\right| = 2\sinh\left|\frac{1}{2}\log(\rho\mathscr{A})\right|,$$

and the $i^{\text{th}}$ eigenvalue of $\log(\rho\mathscr{A})$ is given by

$$-\sum_{k=1}^{m}\frac{\log\lambda_k}{2} + \log\lambda_i.$$

If we choose $k_0$ such that $|\log\lambda_{k_0}| = \max_k\left|\log\lambda_k\right|$, then

$$\left|-\sum_{k=1}^{m}\frac{\log\lambda_k}{2} + \log\lambda_i\right| \leqslant \frac{m}{2}\left|\log\lambda_{k_0}\right|.$$

Hence,

$$\sigma_{\max}\left(\left|\hat{\mathbf{S}}_{g,h}\right|\right) \leqslant 2\sinh\frac{m}{4}\left|\log\lambda_{k_0}\right| = \delta_{g,h}(x),$$

justifying the definition of $\delta_{g,h}$. A similar calculation shows the assertion for $\mathbf{S}_{g,h}$. ∎

The following Lemma provides the trace class operator required in the decomposition formula in assumption (4) of the Belopol'skii-Birman Theorem A.3.

**Lemma 5.3.** *Let $g, h \in \mathrm{Metr}M$, $g \sim h$. We define the unbounded operator*

$$\begin{aligned}
T_s^{g,h} &: \Omega_{\mathsf{L}^2}(M,g) \to \Omega_{\mathsf{L}^2}(M,h)\\
T_s^{g,h} &:= (\hat{P}_s^{h,g})^*\hat{\mathbf{U}}_{g,h}\hat{P}_s^{g} - (\hat{P}_s^{h})^*\mathbf{U}_{g,h}\hat{P}_s^{g,h} + (\check{P}_s^{g,h})^*\hat{\mathbf{U}}_{g,h}\check{P}_s^{g} - (\check{P}_s^{h})^*\mathbf{U}_{g,h}\check{P}_s^{h,g}\\
&\quad - P_s^h\hat{\mathbf{S}}_{g,h;h}\widetilde{\mathbf{U}}_{g,h}\hat{\mathbf{S}}_{g,h;g}P_{s/2}^g\Delta_g P_{s/2}^g.
\end{aligned}$$

*Then we get, for $\alpha_1 \in \operatorname{dom}\Delta_g$, $\alpha_2 \in \operatorname{dom}\Delta_h$ and $s > 0$,*

$$\left\langle\alpha_2, T_s^{g,h}\alpha_1\right\rangle_h = \left\langle\Delta_h\alpha_2, P_s^h I P_s^g\alpha_1\right\rangle_h - \left\langle\alpha_2, P_s^h I P_s^g\Delta_g\alpha_1\right\rangle_h.$$

**Proof.** First note that

$$\mathbf{d}^2 = 0 \qquad \text{and} \qquad \boldsymbol{\delta}^2 = 0. \tag{5.2}$$

Since $\Delta_\nu$ is essentially self-adjoint, for $\nu \in \{g,h\}$, we can assume $\alpha_1 \in \Omega_{\mathsf{C}_c^\infty}(M,g)$ and $\alpha_2 \in \Omega_{\mathsf{C}_c^\infty}(M,h)$ to be compactly supported. Then

$$\begin{aligned}
&\left\langle\Delta_h\alpha_2, P_s^h I P_s^g\alpha_1\right\rangle_h - \left\langle\alpha_2, P_s^h I P_s^g\Delta_g\alpha_1\right\rangle_h\\
&\quad = \left\langle\Delta_h P_s^h\alpha_2, I P_s^g\alpha_1\right\rangle_h - \left\langle\Delta_g I^{-1}P_s^h\alpha_2, P_s^g\alpha_1\right\rangle_h - \left\langle P_s^h\alpha_2, \left(I - (I^{-1})^*\right)P_s^g\Delta_g\alpha_1\right\rangle_h
\end{aligned}$$

$$
\begin{aligned}
&= -\left\langle (\mathbf{d} + \boldsymbol{\delta}_h) P_s^h \alpha_2, (\mathbf{d} + \boldsymbol{\delta}_h) I P_s^g \alpha_1 \right\rangle_h + \left\langle (\mathbf{d} + \boldsymbol{\delta}_g) I^{-1} P_s^h \alpha_2, (\mathbf{d} + \boldsymbol{\delta}_g) P_s^g \alpha_1 \right\rangle_h \\
&\quad - \left\langle P_s^h \alpha_2, \left(I - (I^{-1})^*\right) P_s^g \Delta_g \alpha_1 \right\rangle_h \\
&\overset{(5.2)}{=} -\left\langle \mathbf{d}\, P_s^h \alpha_2, \mathbf{d}\, I P_s^g \alpha_1 \right\rangle_h - \left\langle \boldsymbol{\delta}_h P_s^h \alpha_2, \boldsymbol{\delta}_h I P_s^g \alpha_1 \right\rangle_h \\
&\quad + \left\langle \mathbf{d}\, I^{-1} P_s^h \alpha_2, \mathbf{d}\, P_s^g \alpha_1 \right\rangle_h + \left\langle \boldsymbol{\delta}_g I^{-1} P_s^h \alpha_2, \boldsymbol{\delta}_g P_s^g \alpha_1 \right\rangle_h \qquad (5.3) \\
&\quad - \left\langle P_s^h \alpha_2, \left(I - (I^{-1})^*\right) P_s^g \Delta_g \alpha_1 \right\rangle_h .
\end{aligned}
$$

Let us treat the terms separately. For the last term in (5.3),

$$
\begin{aligned}
&\left\langle P_s^h \alpha_2, \left(I - (I^{-1})^*\right) P_s^g \Delta_g \alpha_1 \right\rangle_h \\
&\quad = \int_M \left\langle P_s^h \alpha_2, \left(\mathscr{A}^{-1/2} - (\rho^{-1} \mathscr{A}^{-1/2}\right) P_s^g \Delta_g \alpha_1 \right\rangle \mathrm{d\,vol}_h \\
&\quad = \int_M \left\langle P_s^h \alpha_2, \left(1 - \rho^{-1}\right) \mathscr{A}^{-1/2} P_s^g \Delta_g \alpha_1 \right\rangle \mathrm{d\,vol}_h \\
&\quad = \int_M \left\langle P_s^h \alpha_2, \hat{\mathbf{S}}_{g,h} \rho^{-1/2} \mathscr{A}^{-1/2} P_s^g \Delta_g \alpha_1 \right\rangle \mathrm{d\,vol}_h \\
&\quad = \int_M \left\langle P_s^h \alpha_2, \left|\hat{\mathbf{S}}_{g,h}\right|^{1/2} (\operatorname{sgn} \hat{\mathbf{S}}_{g,h}) (\rho \mathscr{A})^{-1/2} \left|\hat{\mathbf{S}}_{g,h}\right|^{1/2} P_s^g \Delta_g \alpha_1 \right\rangle \mathrm{d\,vol}_h \\
&\quad = \left\langle \alpha_2, P_s^h \hat{\mathbf{S}}_{g,h;h} \widetilde{\mathbf{U}}_{g,h} \hat{\mathbf{S}}_{g,h;g} P_{s/2}^g \Delta_g P_{s/2}^g \alpha_1 \right\rangle_h .
\end{aligned}
$$

For the first and third term involving the exterior derivative $\mathbf{d}$, we get

$$
\begin{aligned}
&\left\langle \mathbf{d}\, I^{-1} P_s^h \alpha_2, \mathbf{d}\, P_s^g \alpha_1 \right\rangle_h - \left\langle \mathbf{d}\, P_s^h \alpha_2, \mathbf{d}\, I P_s^g \alpha_1 \right\rangle_h \\
&\quad = \left\langle (I \mathbf{d}\, I^{-1}) P_s^h \alpha_2, (I^{-1})^* \mathbf{d}\, P_s^g \alpha_1 \right\rangle_h - \left\langle \mathbf{d}\, P_s^h \alpha_2, I (I^{-1} \mathbf{d}\, I) P_s^g \alpha_1 \right\rangle_h \\
&\quad = \left\langle \hat{P}_s^{h,g} \alpha_2, (I^{-1})^* \mathbf{d}\, P_s^g \alpha_1 \right\rangle_h - \left\langle \mathbf{d}\, P_s^h \alpha_2, I \hat{P}_s^{g,h} \alpha_1 \right\rangle_h \\
&\quad = \left\langle \alpha_2, (\hat{P}_s^{h,g})^* (I^{-1})^* \mathbf{d}\, P_s^g \alpha_1 - (\mathbf{d}\, P_s^h)^* I \hat{P}_s^{g,h} \alpha_1 \right\rangle_h \\
&\quad = \int \left( \alpha_2, (\hat{P}_s^{h,g})^* \rho^{-1} \mathscr{A}^{-1/2} \mathbf{d}\, P_s^g \alpha_1 - (\mathbf{d}\, P_s^h)^* \mathscr{A}^{-1/2} \hat{P}_s^{g,h} \alpha_1 \right) \mathrm{d\,vol}_h \\
&\quad = \left\langle \alpha_2, \left( (\hat{P}_s^{h,g})^* \hat{\mathbf{U}}_{g,h} \mathbf{d}\, P_s^g - (\mathbf{d}\, P_s^h)^* \mathbf{U}_{g,h} \hat{P}_s^{g,h} \right) \alpha_1 \right\rangle_h \\
&\quad = \left\langle \alpha_2, \left( (\hat{P}_s^{h,g})^* \hat{\mathbf{U}}_{g,h} \hat{P}_s^g - (\hat{P}_s^h)^* \mathbf{U}_{g,h} \hat{P}_s^{g,h} \right) \alpha_1 \right\rangle_h .
\end{aligned}
$$

Similarly, for the codifferential $\boldsymbol{\delta}$,

$$
\begin{aligned}
&\left\langle \boldsymbol{\delta}_g I^{-1} P_s^h \alpha_2, \boldsymbol{\delta}_g P_s^g \alpha_1 \right\rangle_h - \left\langle \boldsymbol{\delta}_h P_s^h \alpha_2, \boldsymbol{\delta}_h I P_s^g \alpha_1 \right\rangle_h \\
&\quad = \left\langle (I \boldsymbol{\delta}_g I^{-1}) P_s^h \alpha_2, (I^{-1})^* \boldsymbol{\delta}_g P_s^g \alpha_1 \right\rangle_h - \left\langle \boldsymbol{\delta}_h P_s^h \alpha_2, I (I^{-1} \boldsymbol{\delta}_h I) P_s^g \alpha_1 \right\rangle_h \\
&\quad = \left\langle \check{P}_s^{g,h} \alpha_2, (I^{-1})^* \boldsymbol{\delta}_h P_s^g \alpha_1 \right\rangle_h - \left\langle \boldsymbol{\delta}_h P_s^h \alpha_2, I \check{P}_s^{h,g} \alpha_1 \right\rangle_h \\
&\quad = \left\langle \alpha_2, \left( (\check{P}_s^{g,h})^* \hat{\mathbf{U}}_{g,h} \boldsymbol{\delta}_h P_s^g - (\boldsymbol{\delta}_h P_s^h)^* \mathbf{U}_{g,h} \check{P}_s^{h,g} \right) \alpha_1 \right\rangle_h \\
&\quad = \left\langle \alpha_2, \left( (\check{P}_s^{g,h})^* \hat{\mathbf{U}}_{g,h} \check{P}_s^g - (\check{P}_s^h)^* \mathbf{U}_{g,h} \check{P}_s^{h,g} \right) \alpha_1 \right\rangle_h . \qquad \blacksquare
\end{aligned}
$$

We are finally in the position to proof our main result.

**Proof of Theorem 4.3.** We check the assumptions of the Belopol'skii-Birman theorem A.3.

Since $g \sim h$, the operator $I \equiv I_{g,h}$ is well-defined and bounded and has a bounded inverse $I^{-1} \equiv I_{h,g}$, so (1) follows. By Lemma 5.1, also assumption (2) is satisfied.

Recalling that by (2.14), $I^*_{g,h} = \rho_{g,h} I^{-1}_{g,h}$, we see that the operator $(I^*I - 1)\mathrm{e}^{s\Delta_g}$ has the integral kernel

$$\begin{aligned}\left[(I^*I-1)\mathrm{e}^{s\Delta_g}\right](x,y) &= \left(\rho_{g,h}-1\right) p^g_s(x,y)\\ &= \rho^{1/2}_{g,h}(\operatorname{sgn}\mathbf{S}_{g,h})\left|\mathbf{S}_{g,h}\right|^{1/2}\left|\mathbf{S}_{g,h}\right|^{1/2} p^g_s(x,y).\end{aligned}$$

Thus by Lemma 3.5 again, for some $s > 0$,

$$\begin{aligned}\int\left|\left[(I^*I-1)\mathrm{e}^{s\Delta_g}\right](x,y)\right|^2 \mathrm{vol}_g(\mathrm{d}y) &\leqslant \left\|\rho^{1/2}_{g,h}\mathbf{S}_{g,h}\right\|_\infty \left|\mathbf{S}_{g,h}\right| \int p^g_s(x,y)^2\,\mathrm{vol}_g(\mathrm{d}y)\\ &\leqslant C(\gamma, c_\gamma, s)\left\|\rho^{1/2}_{g,h}\mathbf{S}_{g,h}\right\|_\infty \left|\mathbf{S}_{g,h}\right| \Phi_g(x,s)\int p^{g,(0)}_s(x,y)\,\mathrm{vol}_g(\mathrm{d}y),\end{aligned}$$

and we arrive at the Hilbert-Schmidt estimate

$$\iint\left|\left[(I^*I-1)\mathrm{e}^{s\Delta_g}\right](x,y)\right|^2 \mathrm{vol}_g(\mathrm{d}y)\,\mathrm{vol}_g(\mathrm{d}x) \lesssim \int \delta_{g,h}(x)\Phi_g(x,s)\,\mathrm{vol}_g(\mathrm{d}x) < \infty.$$

So, the operator $(I^*I - 1)\mathrm{e}^{s\Delta_g}$ is Hilbert-Schmidt, hence compact, proving assumption (3).

Finally, we prove (4). Using Lemma 5.3 it remains to show that $T^{g,h}_s$ is trace class. Since the product of Hilbert-Schmidt operators is trace class, we prove that the operators $\hat{P}^\nu_s$, $\check{P}^\nu_s$, for $\nu \in \{g, h\}$, and $\hat{P}^{g,h}_s$, $\check{P}^{g,h}_s$, $\check{P}^{h,g}_s$ are Hilbert-Schmidt. Recall that $\hat{p}^\nu_s(x,y)$, $\check{p}^\nu_s(x,y)$, $\hat{p}^{g,h}_s(x,y)$, $\check{p}^{g,h}_s(x,y)$ and $\check{p}^{h,g}_s(x,y)$ are the corresponding jointly smooth integral kernel of $\hat{P}^\nu_s$, $\check{P}^\nu_s$, $\hat{P}^{g,h}_s$, $\check{P}^{g,h}_s$ and $\check{P}^{h,g}_s$, respectively.

Then, by (4.1),

$$\iint\left|\hat{p}^\nu_s(x,y)\right|^2_\nu \mathrm{vol}_\nu(\mathrm{d}y)\,\mathrm{vol}_\nu(\mathrm{d}x) \leqslant \int \Psi_\nu(x,s)\Phi_\nu(x,s)\,\mathrm{vol}_\nu(\mathrm{d}x)$$

and, by (4.2),

$$\iint\left|\check{p}^\nu_s(x,y)\right|^2_\nu \mathrm{vol}_\nu(\mathrm{d}y)\,\mathrm{vol}_\nu(\mathrm{d}x) \leqslant \int \Psi_\nu(x,s)\Phi_\nu(x,s)\,\mathrm{vol}_\nu(\mathrm{d}x)$$

So similarly, we have by (4.3)

$$\iint\left|\hat{p}^{g,h}_s(x,y)\right|^2_g \mathrm{vol}_g(\mathrm{d}y)\,\mathrm{vol}_g(\mathrm{d}x) \lesssim \int\left(\delta^\nabla_{g,h}(x) + \Xi_g(x,s)\right)\Phi_g(x,s)\,\mathrm{vol}_g(\mathrm{d}x)$$

and by (4.4)

$$\iint\left|\check{p}^{g,h}_s(x,y)\right|^2_g \mathrm{vol}_g(\mathrm{d}y)\,\mathrm{vol}_g(\mathrm{d}x) \lesssim \int\left(\delta^\nabla_{g,h}(x) + \Xi_g(x,s)\right)\Phi_g(x,s)\,\mathrm{vol}_g(\mathrm{d}x)$$

Finally, by (4.5),

$$\iint\left|\check{p}^{h,g}_s(x,y)\right|^2_g \mathrm{vol}_g(\mathrm{d}y)\,\mathrm{vol}_g(\mathrm{d}x) \lesssim \int\left(\delta^\nabla_{g,h}(x) + \Xi_g(x,s)\right)\Phi_g(x,s)\,\mathrm{vol}_g(\mathrm{d}x)$$

This completes the proof. ∎

# 6 Applications and examples

## 6.1 Ricci flow

We first generalise a result, given in [GT20], concerning the stability of the absolutely continuous spectrum of a family of metrics evolving under a Ricci flow. Let therefore $\mathsf{R}_g$ be the Riemannian curvature tensor with respect to the metric $g$.

**Corollary 6.1.** *Let $S > 0$, $\lambda \in \mathbb{R}$ and assume that*

*(a) the family $(g_s)_{0 \leqslant s \leqslant S} \subset \mathrm{Metr} M$ evolves under a Ricci-type flow*

$$\partial_s g_s = \lambda \, \mathrm{Ric}_{g_s}, \qquad \text{for all } \; 0 \leqslant s \leqslant S$$

*(b) the initial metric $g_0$ is geodesically complete*

*(c) there is some $C > 0$ such that $\left|\mathsf{R}_{g_s}\right|_{g_s}, \left|\nabla^{g_s} \mathsf{R}_{g_s}\right|_{g_s} \leqslant C \qquad$ for all $\; 0 \leqslant s \leqslant S$.*

*We set, for all $x \in M$,*

$$\begin{aligned} M_1(x) &:= \sup\left\{ \left|\mathrm{Ric}_{g_s}(v,v)\right| \,:\, 0 \leqslant s \leqslant S,\; v \in \mathsf{T}_x M,\; |v|_{g_s} \leqslant 1 \right\}, \\ M_2(x) &:= \sup\Big\{ \left|\nabla^{g_s}_v \mathrm{Ric}_{g_s}(u,w) + \nabla^{g_s}_u \mathrm{Ric}_{g_s}(v,w) + \nabla^{g_s}_w \mathrm{Ric}_{g_s}(u,v)\right| \,:\, 0 \leqslant s \leqslant S, \\ &\qquad\qquad u,v,w \in \mathsf{T}_x M,\; |u|_{g_s}, |v|_{g_s}, |w|_{g_s} \leqslant 1 \Big\}. \end{aligned}$$

*Let $\mathsf{B}_g(x,R)$ denote the open geodesic ball. If*

$$\int \mathrm{vol}_{g_0}(\mathsf{B}_{g_0}(x,1))^{-1} \max\left\{ \sinh\left(\frac{m}{4} S \,|\lambda|\, M_1(x)\right), M_2(x) \right\} \mathrm{vol}_{g_0}(\mathrm{d}x) < \infty,$$

*then $\sigma_{\mathrm{ac}}(\Delta_{g_s}) = \sigma_{\mathrm{ac}}(\Delta_{g_0})$ for all $0 \leqslant s \leqslant S$.*

**Proof.** The Ricci flow equation together with (a) implies that $g_s \sim g_0$ for all $0 \leqslant s \leqslant S$ and all $g_s$ are complete. Assumption (c) assures that $\Xi(x,s)$ is bounded. By the same arguments as in [GT20, Corollary B],

$$\delta_{g_s, g_0} \leqslant \sinh\left(\frac{m}{4} S \,|\lambda|\, M_1(x)\right)$$

and as in the proof of [BG20, Theorem 6.1]

$$\delta^{\nabla}_{g_s, g_0}(x) \leqslant C M_2(x)$$

and so the claim follows. ■

## 6.2 Differential $k$-forms

We specify our main result, Theorem 4.3, to differential $k$-forms. Set

$$\begin{aligned} \overline{K}^{(k)}_g(x) &:= \max\left\{ \left(\mathscr{R}^{(k)}_g v, v\right) \,:\, v \in \textstyle\bigwedge^k \mathsf{T}_y(M,g), |v|_g = 1, y \in \mathsf{B}_g(x,1) \right\}, \\ \underline{K}^{(k)}_g(x) &:= \min\left\{ \left(\mathscr{R}^{(k)}_g v, v\right) \,:\, v \in \textstyle\bigwedge^k \mathsf{T}_y(M,g), |v|_g = 1, y \in \mathsf{B}_g(x,1) \right\}, \end{aligned}$$

for the corresponding constants defined analogously to (3.10) and (3.11), respectively. Following the same lines as the proof of our main result, Theorem 4.3, we get the following

**Corollary 6.2.** *Let* $g, h \in \mathrm{Metr}M$*,* $g \sim h$*, and assume that there exists* $C < \infty$ *such that* $\left|\delta^{\nabla}_{g,h}\right| \leqslant C$ *and that for both* $\nu \in \{g, h\}$*, we have* $\left|\mathscr{R}_\nu\right|_\nu \in \mathsf{K}(M)$ *and*

$$\int \max\left\{\delta_{g,h}(x), \delta^{\nabla}_{g,h}(x) + \Xi^{(k)}_g(x,s), \Psi^{(k)}_\nu(x,s)\right\} \Phi_\nu(x,s)\,\mathrm{vol}_\nu(\mathrm{d}x) < \infty, \qquad \textit{some } s > 0, \tag{6.1}$$

*where*

$$\Psi^{(k),\pm}_\nu(x,s) := \Psi^{(k),+}_\nu(x,s) \wedge \Psi^{(k),-}_\nu(x,s),$$
$$\Xi^{(k),\pm}_\nu(x,s) := \Xi^{(k),+}_\nu(x,s) \wedge \Xi^{(k),-}_\nu(x,s),$$

*with*

$$\Psi^{(k),\pm}_\nu(x,s) := \mathrm{e}^{D(\gamma, c_\gamma(\underline{\mathscr{R}}^-), c_q^{1/q})s \,+\, (\overline{K}^{(k)}_\nu(x) + \underline{K}^{(k\pm1)}_\nu(x))^- s/2}\, \mathrm{e}^{C(m, \underline{K}^{(0),-}_\nu(x))s/2} s^{-1}.$$

*and*

$$\Xi^{(k),\pm}_\nu(x,s) := \Psi^{(k),\pm}_\nu(x,s) + s^{3/2}\Psi^{(k),\pm}_\nu(x,s) \max_{y \in \mathsf{B}(x,1)} |\nabla \mathsf{R}(y)|\,.$$

*For all* $0 \leqslant k \leqslant m$*, let*

$$I^{(k)} := I^{(k)}_{g,h} : \Omega^k_{\mathsf{L}^2}(M, g) \to \Omega^k_{\mathsf{L}^2}(M, h), \quad \alpha \mapsto \bigwedge\nolimits^k A^{-1/2}(\alpha)$$

*be the bounded identification operator acting on* $k$*-forms. Then, for all* $0 \leqslant k \leqslant m$*, the wave operators*

$$W_\pm(\Delta^{(k)}_h, \Delta^{(k)}_g, I^{(k)}) = \operatorname*{s\text{-}lim}_{t \to \pm\infty} \mathrm{e}^{it\Delta^{(k)}_h} I \mathrm{e}^{-it\Delta^{(k)}_g} \mathsf{P}_{\mathrm{ac}}(\Delta^{(k)}_g)$$

*exist and are complete. Moreover,* $W_\pm(\Delta^{(k)}_h, \Delta^{(k)}_g, I^{(k)})$ *are partial isometries with initial space* $\operatorname{ran}\mathsf{P}_{\mathrm{ac}}(\Delta^{(k)}_g)$ *and final space* $\operatorname{ran}\mathsf{P}_{\mathrm{ac}}(\Delta^{(k)}_h)$*, and we have* $\sigma_{\mathrm{ac}}(\Delta^{(k)}_g) = \sigma_{\mathrm{ac}}(\Delta^{(k)}_h)$*.*

**Proof.** We omit the metric in the notation. By similar calculations as in the proofs of Theorem 3.9 and Theorem 3.10, we see that for every $\alpha \in \Omega^k_{\mathsf{L}^2}(M)$,

$$\left|(\mathbf{d}^{(k)} P_s \alpha)_x\right|^2 \leqslant \Psi^{(k),+}(x,s)\Phi(x,s)\,\|\alpha\|^2_{\Omega^k_{\mathsf{L}^2}(M)},$$
$$\left|(\boldsymbol{\delta}^{(k)} P_s \alpha)_x\right|^2 \leqslant \Psi^{(k),-}(x,s)\Phi(x,s)\,\|\alpha\|^2_{\Omega^k_{\mathsf{L}^2}(M)},$$
$$\left|\left(\nabla P_s \alpha, \xi\right)\right|^2 \leqslant \max\left\{\Xi^{(k),+}(x,s), \Xi^{(k),-}(x,s)\right\} \Phi(x,s)\,\|\alpha\|^2_{\Omega^k_{\mathsf{L}^2}(M)},$$

Noticing that

$$\Delta_\nu = \bigoplus_{k=0}^m \Delta^{(k)}_\nu \qquad \text{and} \qquad I = \bigoplus_{k=0}^m I^{(k)}$$

the proof now follows the lines of the proof of our main result, Theorem 4.3. ∎

## 6.3 Conformal perturbations

We study the important case of conformally equivalent metrics: Given a smooth function $\psi : M \to \mathbb{R}$, we define peturbated metric by $g_\psi := \mathrm{e}^{2\psi} g$. Note that $g$ and $g_\psi$ are quasi-isometric, if and only if $\psi$ is bounded (cf. Example 2.7 above).

The bounded identification operator is now given by

$$I := I_{g,g_\psi} : \Omega_{\mathsf{L}^2}(M,g) \to \Omega_{\mathsf{L}^2}(M,g_\psi)$$
$$I_{g,g_\psi}\eta(x) \mapsto \mathrm{e}^{-\psi(x)}\eta(x).$$

Given a smooth function $\psi$ on $M$, we define

$$\tau := \bigoplus_{k=0}^{m}(m-2k)\mathbf{1}_{\bigwedge^k \mathsf{T}^* M} \in \mathscr{D}^{(0)}(M; \textstyle\bigwedge \mathsf{T}^* M),$$
$$\mathrm{e}^{\psi\tau} := \bigoplus_{k=0}^{m}\mathrm{e}^{(m-2k)\psi}\mathbf{1}_{\bigwedge^k \mathsf{T}^* M} \in \mathscr{D}^{(0)}(M; \textstyle\bigwedge \mathsf{T}^* M).$$

Next, we collect some useful transformation rules for the conformal metric $g_\psi$ in terms of $g$. A standard reference for various invariants of conformal metric change in part (a) is [Bes87, 1.159 Theorem].

**Proposition 6.3.** *Let $\psi : M \to \mathbb{R}$ be smooth.*

*(a) We have*

$$(\cdot,\cdot)^{(k)}_{g_\psi} = \mathrm{e}^{-2k\psi}(\cdot,\cdot)^{(k)}_g \qquad \text{for all } k \in \{0,\dots,m\} \tag{6.2a}$$

$$\mathrm{d}\,\mathrm{vol}_{g_\psi} = \mathrm{e}^{m\psi}\mathrm{d}\,\mathrm{vol}_g \tag{6.2b}$$

$$\bullet \lrcorner_{g_\psi} \alpha = \mathrm{e}^{-2\psi}\left(\bullet \lrcorner_g \alpha\right) \qquad \text{for all } \alpha \in \Omega^1(M) \tag{6.2c}$$

$$\nabla^{g_\psi}_X Y = \nabla^g_X Y + \mathbf{d}\psi(X)Y + \mathbf{d}\psi(Y)X - (X,Y)_g \operatorname{grad}_g \psi \quad \text{for all } X,Y \in \Gamma_{\mathsf{C}^\infty}(\mathsf{T}M) \tag{6.2d}$$

$$\boldsymbol{\delta}_{g_\psi}\alpha = \mathrm{e}^{-2\psi}(\boldsymbol{\delta}_g\alpha - \tau\,(\mathbf{d}\psi)^{\sharp_g} \lrcorner_g \alpha) \qquad \text{for all } \alpha \in \Omega^k(M) \tag{6.2e}$$

*(b) If $\psi$ is bounded, then*

$$I^* = \mathrm{e}^{m\psi}I^{-1}. \tag{6.3}$$

**Remark 6.4.** We note that the canonical musical isomorphisms $\sharp$ and $\flat$ between $\mathsf{T}M$ and $\mathsf{T}^*M$ do not agree for $g$ and $g_\psi$.

**Proof.** (a) We show (6.2e). Using (6.2a) and (6.2b), for any $\eta_1 \in \Omega^{k-1}(M)$, $\eta_2 \in \Omega^k(M)$,

$$\begin{aligned}
\left\langle \eta_1, \boldsymbol{\delta}_{g_\psi}\eta_2 \right\rangle_{g_\psi} &= \left\langle \eta_1, \mathrm{e}^{-2\psi}\left(\boldsymbol{\delta}_g\eta_2 + (m-2p)(\mathbf{d}\psi)^{\sharp_g} \lrcorner\, \eta_2\right)\right\rangle_{g_\psi} \\
&= \left\langle \mathrm{e}^{(m-2(p-1))\psi}\eta_1, \mathrm{e}^{-2\psi}\boldsymbol{\delta}_g\eta_2 \right\rangle_g + \left\langle \mathrm{e}^{(m-2(p-1))\psi}\eta_1, \mathrm{e}^{-2\psi}(m-2p)(\mathbf{d}\psi)^{\sharp_g} \lrcorner\, \eta_2 \right\rangle_g \\
&= \left\langle \mathbf{d}\left(\mathrm{e}^{(m-2p)\psi}\eta_1\right), \eta_2 \right\rangle_g + \left\langle (m-2p)\mathrm{e}^{(m-2p)\psi}\eta_1, (\mathbf{d}\psi)^{\sharp_g} \lrcorner\, \eta_2 \right\rangle_g \\
&= \left\langle \mathrm{e}^{(m-2p)\psi}\mathbf{d}\eta_1, \eta_2 \right\rangle_g + \left\langle \mathbf{d}\left(\mathrm{e}^{(m-2p)\psi}\right) \wedge \eta_1, \eta_2 \right\rangle_g - \left\langle (m-2p)\mathrm{e}^{(m-2p)\psi}\mathbf{d}\psi \wedge \eta_1, \eta_2 \right\rangle_g \\
&= \langle \mathbf{d}\eta_1, \eta_2 \rangle_{g_\psi}\ .
\end{aligned}$$

(b) Follows from (2.14). ■

**Theorem 6.5.** *Let $\psi : M \to \mathbb{R}$ be smooth with $\psi$ bounded, and assume that $g, g_\psi \in \mathrm{Metr}M$ with $g_\psi = \mathrm{e}^{2\psi} g$ such that $\left|\delta^{\nabla}_{g,g_\psi}\right| \leqslant C$ for some $C < \infty$ and that for both $v \in \{g, g_\psi\}$, we have $\left|\mathscr{R}_v\right|_v \in \mathrm{K}(M)$ and*

$$\int \max\left\{\sinh\left|\frac{m}{4}\psi(x)\right|, \delta^{\nabla}_{g,h}(x) + \Xi_g(x,s), \Psi_v(x,s)\right\}\Phi_v(x,s)\,\mathrm{vol}_v(\mathrm{d}x) < \infty, \qquad \text{some } s > 0. \tag{6.4}$$

*Then the wave operators*

$$W_\pm(\Delta_{g_\psi}, \Delta_g, I) = \operatorname*{s-lim}_{t\to\pm\infty} \mathrm{e}^{it\Delta_{g_\psi}} I \mathrm{e}^{-it\Delta_g} \mathrm{P}_{\mathrm{ac}}(\Delta_g)$$

*exist and are complete. Moreover, $W_\pm(\Delta_{g_\psi}, \Delta_g, I)$ are partial isometries with initial space* ran $\mathrm{P}_{\mathrm{ac}}(\Delta_g)$ *and final space* ran $\mathrm{P}_{\mathrm{ac}}(\Delta_{g_\psi})$*, and we have $\sigma_{\mathrm{ac}}(\Delta_g) = \sigma_{\mathrm{ac}}(\Delta_g)$.*

**Proof.** Using Example 2.7, we have $\delta_{g,g_\psi} = 2\sinh\frac{m}{4}|\psi|$ and

$$g \sim g_\psi \qquad \Longleftrightarrow \qquad \psi \text{ bounded.}$$

Hence the claim follows from our main result, Theorem 4.3. ■

**Remark 6.6.** By the same argument as in the Proof of Corollary 6.2, we can deduce Theorem 6.5 for the wave operators acting on $k$-forms but with appropriate localised constants respecting the degree of the differential form (cf. Proof of Corollary 6.2 above).

## 6.4 Global curvature bounds

Let $\mathrm{R}_v$ be the Riemannian curvature tensor with respect to the metric $v \in \{g, h\}$. Then the curvature operator

$$Q_v \in \mathscr{D}^{(0)}\left(M; \bigwedge^2 \mathrm{T}^*M\right)$$

is self-adjoint and uniquely determined by the equation

$$\left(Q_v(X\wedge Y), U\wedge V\right)_v = \left(\mathrm{R}_v(X,Y)U, V\right)_v$$

for all smooth vector fields $X, Y, U, V \in \Gamma_{C^\infty}(\mathrm{T}M)$.

By the Gallot-Meyer estimate [GM75], a global bound $Q_v \geqslant -K$, for some constant $K > 0$, already implies that curvature endomorphism in the Weitzenböck formula (3.1) is globally bounded by

$$\mathscr{R}^{(k)}_v \geqslant -Kk(m-k). \tag{6.5}$$

**Remark 6.7.** In particular, if $Q_v \geqslant -K$, for some constant $K > 0$, then for $k = 1$, we have

$$\mathrm{Ric}_v \geqslant -K(m-1).$$

We set, for $v \in \{g,h\}$,

$$\Theta_v(x,s) := \left(1 + \max_{y\in \mathsf{B}_v(x,1)} \left|\nabla^v \mathsf{R}_v(y)\right|\right)^2 .$$

We get the following consequential result.

**Theorem 6.8.** *Let $Q_v \geqslant -K$, for some constant $K > 0$ for both $v \in \{g,h\}$. Let $g,h \in \mathrm{Metr}M$ such that $g \sim h$ and assume that there exists $C < \infty$ such that $\left|\delta^{\nabla}_{g,h}\right| \leqslant C$ and that for some (then both by quasi-isometry) $v \in \{g,h\}$*

$$\int \max\left\{\delta_{g,h}(x), \delta^{\nabla}_{g,h}(x) + \Theta_g(x,s)\right\} \Phi_v(x,s)\,\mathrm{vol}_v(\mathrm{d}x) < \infty, \qquad \text{some } s > 0.$$

*Then, the wave operators $W_\pm(\Delta_h, \Delta_g, I)$ exist and are complete. Moreover, $W_\pm(\Delta_h, \Delta_g, I)$ are partial isometries with initial space $\operatorname{ran} \mathsf{P}_{\mathrm{ac}}(\Delta_g)$ and final space $\operatorname{ran} \mathsf{P}_{\mathrm{ac}}(\Delta_h)$, and we have $\sigma_{\mathrm{ac}}(\Delta_g) = \sigma_{\mathrm{ac}}(\Delta_h)$.*

**Remark 6.9.** Following the arguments in the previous subsections 6.4 and 6.3, we may elaborate Corollary 6.2 and Theorem 6.5 for the special case of differential $k$-forms and the conformal metric change, respectively, in the same manner.

**Lemma 6.10.** *Let $g,h \in \mathrm{Metr}M$, $g \sim h$ with $Q_v \geqslant -K$, for some constant $K > 0$ for both $v \in \{g,h\}$ and assume that the function $\delta^{\nabla}_{g,h}$ is bounded. Then*

$$I_{g,h} \operatorname{dom} \mathbf{q}_g = \operatorname{dom} \mathbf{q}_h.$$

**Proof.** Note that, for any $v \in \{g,h\}$, $\operatorname{dom} \mathbf{q}_v$ is the closure of compactly supported forms $\Omega_{\mathsf{C}_c^\infty}(M,v)$ with respect to the Dirac graph norm

$$\alpha \mapsto \left(\|\alpha\|_v^2 + \left\|\mathbf{D}_v \alpha\right\|_v^2\right)^{1/2} .$$

Let $D$ be a positive constant whose value might change from line to line. For all compactly supported $\alpha \in \Omega_{\mathsf{C}_c^\infty}(M,g)$, we get by the Weitzenböck formula $\Delta_v = \square_v - \mathscr{R}_v$,

$$\begin{aligned} \left\|\mathbf{D}_v I\alpha\right\|_v^2 &= \left\langle \mathbf{D}_v^2 I\alpha, \alpha\right\rangle_v = \left\langle\left(-(\nabla^h)^*\nabla^h - \mathscr{R}_v\right) I\alpha, \alpha\right\rangle_v \\ &= \int \left|\nabla^h (I\alpha)\right|_h^2 \mathrm{d}\,\mathrm{vol}_v - \int \left(\mathscr{R}_v I\alpha, \alpha\right)_v \mathrm{d}\,\mathrm{vol}_v . \end{aligned}$$

By assumption the metrics are quasi-isometric and their Weitzenböck curvature term is bounded, so the second term is bounded by $D\,\|\alpha\|_v^2$.

Following the lines of the proof of Lemma 5.1, we find

$$\left\|\nabla^h (I\alpha)\right\|_h^2 \leqslant D\left(\|\alpha\|_g^2 + \|\nabla^g \alpha\|_g^2\right).$$

Using the Weitzenböck formula once more, we get

$$\begin{aligned} \|\nabla^g \alpha\|_g^2 = \langle \nabla^g\alpha, \nabla^g\alpha\rangle_g &= \left\langle\left(-(\nabla^g)^*\nabla^g - \mathscr{R}\right)\alpha, \alpha\right\rangle_g + \langle \mathscr{R}\alpha, \alpha\rangle_g \\ &\leqslant \left\langle \mathbf{D}_g^2 \alpha, \alpha\right\rangle_g + K\,\langle\alpha,\alpha\rangle_g \\ &\leqslant D\left(\left\|\mathbf{D}_g\alpha\right\|_g^2 + \|\alpha\|_g^2\right). \end{aligned}$$

Hence, we arrive at the estimate

$$\|I\alpha\|_h^2 + \left\|\mathsf{D}_h I\alpha\right\|_h^2 \leqslant D\left(\|\alpha\|_g^2 + \left\|\mathsf{D}_g\alpha\right\|_g^2\right),$$

proving

$$I \operatorname{dom} \mathbf{q}_g \subset \operatorname{dom} \mathbf{q}_h.$$

Since $I^{-1} = I_{g,h}^{-1} = I_{h,g}$ and the arguments above are symmetric in $g$ and $h$, this shows the claim. ■

**Proof of Theorem 6.8.** We omit the metric in the notation. By assumption $\mathscr{R}$ is bounded from below, so the tensor $\widetilde{\mathscr{R}}$, Ric and R are also bounded. In particular $(M, g)$ is stochastically complete, i.e. $\zeta(x) = \infty$ $\mathbb{P}$-a.s. By Gronwall's inequality, we have $\left|\mathscr{Q}_s\right|_{\mathrm{op}}, \left|\mathscr{Q}_s^{-1}\right|_{\mathrm{op}}, \left|\widetilde{\mathscr{Q}}_s\right|_{\mathrm{op}} \leqslant \mathrm{e}^{K^- s/2}$. Following the lines of the proof of Theorem 3.10, we get by Cauchy-Schwarz

$$\begin{aligned}\left|\left(\nabla P_s\alpha(x), \xi\right)\right| &= \left|\mathbb{E}\left(\mathscr{Q}_s^{\mathrm{tr}} /\!/_s^{-1}\alpha(X_s(x)), U_{s\wedge\tau}^{\ell}\right)\right| \\ &\leqslant C(K^-, s)\left[\mathbb{E}\left|\alpha(X_s(x))\right|^2\right]^{1/2}\left(\left[\mathbb{E}\left(\ell_{s\wedge\tau}^{(1)}\right)^2\right]^{1/2} + \left[\mathbb{E}\left(\ell_{s\wedge\tau}^{(2)}\right)^2\right]^{1/2}\right).\end{aligned}$$

By (3.16) and (3.17) the first summand in the bracket is bounded by $C(K^-, s)\,|\xi|$. By (3.21) and (3.22) the second summand is bounded by $C(K^-, s)\max_{y\in \mathsf{B}(x,1)}|\nabla \mathsf{R}(y)|\,|\xi|$. Hence,

$$\left|\left(\nabla P_s\alpha(x), \xi\right)\right|^2 \leqslant C(K^-, s)\,|\xi|^2\,\Phi(x, s)\left(1 + \max_{y\in\mathsf{B}(x,1)}|\nabla \mathsf{R}(y)|\right)^2 \|\alpha\|_{\mathsf{L}^2(M)}^2,$$

and analogously

$$\begin{aligned}\left|(\mathbf{d}P_s\alpha)_x\right|^2 &\leqslant C(K^-, s)\Phi(x, s)\,\|\alpha\|_{\Omega_{\mathsf{L}^2}(M)}^2, \\ \left|(\boldsymbol{\delta}P_s\alpha)_x\right|^2 &\leqslant C(K^-, s)\Phi(x, s)\,\|\alpha\|_{\Omega_{\mathsf{L}^2}(M)}^2.\end{aligned}$$

Following the lines of the proof of Theorem 4.3 we see that the assumptions of Belopol'skii-Birman theorem A.3 are satisfied except making use of Lemma 6.10 (instead of Lemma 5.1) for assumption (2). ■

## 6.5 $\varepsilon$-close Riemannian metrics

In this section, we denote by $\boldsymbol{\kappa}_g$ the sectional curvature with respect to a smooth, complete Riemannian metric $g$.

In [CFG92, Theorem 1.3 & 1.7], Cheeger, Fukaya and Gromov show what is also known as Cheeger-Gromov's thick/thin decomposition.

**Theorem 6.11.** *For all $\varepsilon > 0$ and $n \in \mathbb{Z}_+$, there exists $\xi > 0$ and $k \in \mathbb{Z}_+$ such that if $(M, g)$ is a complete Riemannian manifold with $\left|\boldsymbol{\kappa}_g\right| \leqslant 1$, then there is a $(\xi, k)$-round metric, $g_\varepsilon$, on $M$, such that*

*(i) the Riemannian metric $g_\varepsilon$ is $\varepsilon$-quasi-isometric to $g$, i.e. $(1/C^\varepsilon)g_\varepsilon \leqslant g \leqslant C^\varepsilon g_\varepsilon$*

*(ii) it has bounded covariant derivatives $|\nabla^{g_\varepsilon} - \nabla^g| < \varepsilon$*

*(iii)* $|\left(\nabla^{g_\varepsilon}\right)^k \mathrm{R}_{g_\varepsilon}| < C(m,k,\varepsilon)$*, where the constant* $C$ *depends in addition on the order of derivative* $k$ *and* $\varepsilon$*.*

Assuming that the sectional curvature $\boldsymbol{\kappa}_g$ is bounded by 1, implies that the Riemannian curvature tensor $\mathrm{R}_g$ is bounded, and hence, the curvature operator $Q_g$. Following our results in Section 6.4, we may get

**Theorem 6.12.** *Let* $|\boldsymbol{\kappa}_g| \leqslant 1$*. Then there exists a Riemannian metric* $g_\varepsilon$ *as in Theorem 6.11 that is* $\varepsilon$*-quasi-isometric metric to* $g$*. If for some* $\nu \in \{g, g_\varepsilon\}$

$$\int \delta_{g,g_\varepsilon}(x)\Phi_\nu(x,s)\,\mathrm{vol}_\nu(\mathrm{d}x) < \infty, \qquad \text{some } s > 0,$$

*then the wave operators* $W_\pm(\Delta_{g_\varepsilon}, \Delta_g, I)$ *exist and are complete. Moreover,* $W_\pm(\Delta_{g_\varepsilon}, \Delta_g, I)$ *are partial isometries with initial space* $\operatorname{ran}\mathrm{P}_{\mathrm{ac}}(\Delta_g)$ *and final space* $\operatorname{ran}\mathrm{P}_{\mathrm{ac}}(\Delta_{g_\varepsilon})$*, and we have* $\sigma_{\mathrm{ac}}(\Delta_g) = \sigma_{\mathrm{ac}}(\Delta_{g_\varepsilon})$*.*

**Proof.** By the previous Theorem 6.11 (i), the assumption $|\boldsymbol{\kappa}_g| \leqslant 1$ assures that for any $\varepsilon > 0$ there exists a Riemannian metric $g_\varepsilon$ that is $\varepsilon$-quasi-isometric metric to $g$. Hence,

$$\sup \delta_{g,g_\varepsilon}(x) < \infty \quad \Longleftrightarrow \quad g \sim g_\varepsilon.$$

By Theorem 6.11 (ii), the covariant derivatives are bounded so that

$$\delta^{\nabla}_{g,g_\varepsilon} = |\nabla^{g_\varepsilon} - \nabla^g|_g^2 < \varepsilon. \qquad \blacksquare$$

# A Belopol'skii-Birman theorem

We will use a variant of the Belopol'skii-Birman theorem [RS79; Wei80] which is adapted to our special case of two Hilbert space scattering theory originally to be found in [GT20].

Given a self-adjoint operator $H$ in a Hilbert space $\mathscr{H}$ with its operator valued spectral measure $E_H$, one defines the $H$-**absolutely continuous subspace** $\mathscr{H}_{\mathrm{ac}}(H)$ **of** $\mathscr{H}$ to be the space of all $f \in \mathscr{H}$ such that the Borel measure $\|E_H(\cdot)f\|^2$ on $\mathbb{R}$ is absolutely continuous with respect to the Lebesgue measure. Then $\mathscr{H}_{\mathrm{ac}}(H)$ becomes a closed subspace of $\mathscr{H}$ and the restriction $H_{\mathrm{ac}}$ of $H$ to $\mathscr{H}_{\mathrm{ac}}(H)$ is a well-defined self-adjoint operator. The **absolutely continuous spectrum** $\sigma_{\mathrm{ac}}(H)$ **of** $H$ is defined to be the spectrum of $H_{\mathrm{ac}}$.

**Definition A.1.** Let $H_1$ and $H_2$ be self-adjoint operators on Hilbert spaces $\mathscr{H}_1$ and $\mathscr{H}_2$, respectively, $I : \mathscr{H}_1 \to \mathscr{H}_2$ be a bounded operator and $\mathrm{P}_{\mathrm{ac}}(H_1)$ be the projection onto the absolutely continuous subspace of $H_1$. The **(generalised) wave operators** $W_\pm(H_2, H_1, I)$ exist if the strong limits

$$W_\pm(H_2, H_1, I) = \operatorname*{s-lim}_{t\to\pm\infty} \mathrm{e}^{itH_2} I \mathrm{e}^{-itH_1}\mathrm{P}_{\mathrm{ac}}(H_1)$$

exist.

Note that $W_\pm(H_2, H_1, I)$ may not be isometries. A further assumption is that each scattering state looks asymptotically free which is reflected in the next definition.

**Definition A.2.** Suppose that $W_\pm(H_2, H_1, I)$ exist. We say that they are **complete** if and only if

$$\big(\ker W_\pm(H_2, H_1, I)\big)^\perp = \operatorname{ran} \mathsf{P}_{\mathrm{ac}}(H_1), \qquad \overline{\operatorname{ran} W_\pm(H_2, H_1, I)} = \operatorname{ran} \mathsf{P}_{\mathrm{ac}}(H_2).$$

**Theorem A.3** (Belopol'skii-Birman). *For $k = 1, 2$, let $H_k \geqslant 0$ be self-adjoint operators in a Hilbert space $\mathscr{H}_k$ and $\mathsf{P}_{\mathrm{ac}}(H_k)$ the projection onto the absolutely continuous subspace of $\mathscr{H}_k$ corresponding to $H_k$. Assume that $I \in \mathscr{L}(\mathscr{H}_1, \mathscr{H}_2)$ is a bounded operator such that the following assumptions hold:*

*(1) $I$ has a two-sided bounded inverse*

*(2) We have either $I \operatorname{dom} \sqrt{H_1} = \operatorname{dom} \sqrt{H_2}$ or $I \operatorname{dom} H_1 = \operatorname{dom} H_2$*

*(3) The operator $(I^*I - 1)\mathrm{e}^{-sH_1} : \mathscr{H}_1 \to \mathscr{H}_1$ is compact for some $s > 0$*

*(4) There is a trace class operator $T : \mathscr{H}_1 \to \mathscr{H}_2$ and a number $s > 0$ such that for all $\alpha_1 \in \operatorname{dom} H_1$, $\alpha_2 \in \operatorname{dom} H_2$ we have*

$$\langle \alpha_2, T\alpha_1 \rangle_{\mathscr{H}_2} = \left\langle H_2\alpha_2, \mathrm{e}^{-sH_2} I \mathrm{e}^{-sH_1} \alpha_1 \right\rangle_{\mathscr{H}_2} - \left\langle \alpha_2, \mathrm{e}^{-sH_2} I \mathrm{e}^{-sH_1} H_1 \alpha_1 \right\rangle_{\mathscr{H}_2}.$$

*Then the wave operators*

$$W_\pm(H_2, H_1, I) = \operatorname*{s-lim}_{t \to \pm\infty} \mathrm{e}^{itH_2} I \mathrm{e}^{-itH_1} \mathsf{P}_{\mathrm{ac}}(H_1)$$

*exist and are complete. Moreover, $W_\pm(H_2, H_1, I)$ are partial isometries with initial space $\operatorname{ran} \mathsf{P}_{\mathrm{ac}}(H_1)$ and final space $\operatorname{ran} \mathsf{P}_{\mathrm{ac}}(H_2)$.*

**Proof.** In view of Theorem XI.13 from [RS79] and its proof, it remains to show that for every bounded interval $\mathbb{I}$ the operator $(I^*I - 1)E_1(\mathbb{I})$ is compact, and that there exists a trace class operator $D \in \mathscr{J}^1(\mathscr{H}_1, \mathscr{H}_2)$ such that for every bounded interval $\mathbb{I}$ and all $\alpha_1, \alpha_2$ as above we have

$$\langle \varphi, D\alpha_1 \rangle_{\mathscr{H}_2} = \langle H_2\alpha_2, E_2(\mathbb{I}) I E_1(\mathbb{I})\alpha_1 \rangle_{\mathscr{H}_2} - \langle \alpha_2, E_2(\mathbb{I}) I E_1(\mathbb{I}) H_1\alpha_1 \rangle_{\mathscr{H}_2}.$$

However, using that for all self-adjoint operators $A$ and all Borel functions $\varphi, \varphi' : \mathbb{R} \to \mathbb{C}$ we have

$$\varphi(A)\varphi'(A) \subset (\varphi \cdot \varphi')(A), \quad \operatorname{dom}(\varphi(A)\varphi'(A)) = \operatorname{dom}(\varphi(A)\varphi'(A)) \cap \operatorname{dom} \varphi'(A),$$

the required compactness becomes obvious, and furthermore it is easily justified that

$$D := \mathrm{e}^{sH_2} E_2(\mathbb{I}) I \mathrm{e}^{sH_1} E_1(\mathbb{I})$$

has the required trace class property. ■

By a classical result of Kato (cf. [Kat95, X. Perturbation of continuous spectra and unitary equivalence] and [Kat95, p. 534, Theorem 3.5.] therein), Theorem A.3 implies:

**Theorem A.4.** *Assume in the above situation that the wave operators $W_\pm(H_2, H_1, I)$ exist and are complete. Then the operators $H_{1,\mathrm{ac}}$ and $H_{2,\mathrm{ac}}$ are unitarily equivalent. In particular, we have*

$$\sigma_{\mathrm{ac}}(H_1) = \sigma_{\mathrm{ac}}(H_2).$$

*Acknowledgements.* The author would like to express his thanks to Anton Thalmaier and James Thompson for carefully reading the manuscript and stimulating discussions on the subject, and Batu Güneysu and Sebastian Boldt for the helpful discussions improving this paper. Lastly, the anonymous referees for their detailed and stimulating reports largely improving the presentation of the paper.

## References

[AT10] M. Arnaudon and A. Thalmaier. «Li-Yau type gradient estimates and Harnack inequalities by stochastic analysis». In: *Probabilistic approach to geometry*. Vol. 57. Adv. Stud. Pure Math. Math. Soc. Japan, Tokyo, 2010, pp. 29–48.

[Bes87] A. L. Besse. *Einstein manifolds*. Vol. 10. Ergebnisse der Mathematik und ihrer Grenzgebiete (3) [Results in Mathematics and Related Areas (3)]. Springer-Verlag, Berlin, 1987, pp. xii+510.

[BG20] S. Boldt and B. Güneysu. «Scattering Theory and Spectral Stability under a Ricci Flow for Dirac Operators». In: *arXiv* (Mar. 2020), pp. 1–27. arXiv: 2003.10204.

[BGM17] F. Bei, B. Güneysu and J. Müller. «Scattering theory of the Hodge-Laplacian under a conformal perturbation». In: *J. Spectr. Theory* 7.1 (2017), pp. 235–267.

[CFG92] J. Cheeger, K. Fukaya and M. Gromov. «Nilpotent structures and invariant metrics on collapsed manifolds». In: *J. Amer. Math. Soc.* 5.2 (1992), pp. 327–372.

[CTT18] L.-J. Cheng, A. Thalmaier and J. Thompson. «Quantitative $C^1$-estimates by Bismut formulae». In: *J. Math. Anal. Appl.* 465.2 (2018), pp. 803–813.

[DT01] B. K. Driver and A. Thalmaier. «Heat equation derivative formulas for vector bundles». In: *J. Funct. Anal.* 183.1 (2001), pp. 42–108.

[GM75] S. Gallot and D. Meyer. «Opérateur de courbure et laplacien des formes différentielles d'une variété riemannienne». In: *J. Math. Pures Appl. (9)* 54.3 (1975), pp. 259–284.

[GT20] B. Güneysu and A. Thalmaier. «Scattering theory without injectivity radius assumptions, and spectral stability for the Ricci flow». In: *Annales de l'Institut Fourier* 70.1 (2020), pp. 437–456.

[Gün12] B. Güneysu. «On generalized Schrödinger semigroups». In: *J. Funct. Anal.* 262.11 (2012), pp. 4639–4674.

[Gün17] B. Güneysu. *Covariant Schrödinger semigroups on Riemannian manifolds*. Vol. 264. Operator Theory: Advances and Applications. Birkhäuser/Springer, Cham, 2017, pp. xviii+239.

[HPW14] R. Hempel, O. Post and R. Weder. «On open scattering channels for manifolds with ends». In: *J. Funct. Anal.* 266.9 (2014), pp. 5526–5583.

[Hsu02] E. P. Hsu. *Stochastic analysis on manifolds*. Vol. 38. Graduate Studies in Mathematics. American Mathematical Society, Providence, RI, 2002, pp. xiv+281.

[HT94] W. Hackenbroch and A. Thalmaier. *Stochastische Analysis.* Mathematische Leitfäden. [Mathematical Textbooks]. Eine Einführung in die Theorie der stetigen Semimartingale. [An introduction to the theory of continuous semimartingales]. B. G. Teubner, Stuttgart, 1994, p. 560.

[IS13] K. Itô and E. Skibsted. «Scattering theory for Riemannian Laplacians». In: *J. Funct. Anal.* 264.8 (2013), pp. 1929–1974.

[Jos17] J. Jost. *Riemannian geometry and geometric analysis.* Seventh. Universitext. Springer, Cham, 2017, pp. xiv+697.

[Kat95] T. Kato. *Perturbation theory for linear operators.* Classics in Mathematics. Reprint of the 1980 edition. Springer-Verlag, Berlin, 1995, pp. xxii+619.

[Ken87] W. S. Kendall. «The radial part of Brownian motion on a manifold: a semimartingale property». In: *The Annals of Probability* (1987), pp. 1491–1500.

[Mir75] L. Mirsky. «A trace inequality of John von Neumann». In: *Monatsh. Math.* 79.4 (1975), pp. 303–306.

[MS07] W. Müller and G. Salomonsen. «Scattering theory for the Laplacian on manifolds with bounded curvature». In: *J. Funct. Anal.* 253.1 (2007), pp. 158–206.

[RS79] M. Reed and B. Simon. *Methods of modern mathematical physics. III.* Scattering theory. Academic Press [Harcourt Brace Jovanovich, Publishers], New York-London, 1979, pp. xv+463.

[Sim82] B. Simon. «Schrödinger semigroups». In: *Bull. Amer. Math. Soc. (N.S.)* 7.3 (1982), pp. 447–526.

[Str83] R. S. Strichartz. «Analysis of the Laplacian on the complete Riemannian manifold». In: *J. Functional Analysis* 52.1 (1983), pp. 48–79.

[Tha16] A. Thalmaier. «Geometry of subelliptic diffusions». In: *Geometry, analysis and dynamics on sub-Riemannian manifolds. Vol. II.* EMS Ser. Lect. Math. Eur. Math. Soc., Zürich, 2016, pp. 85–169.

[Tha97] A. Thalmaier. «On the differentiation of heat semigroups and Poisson integrals». In: *Stochastics Stochastics Rep.* 61.3-4 (1997), pp. 297–321.

[TW04] A. Thalmaier and F.-Y. Wang. «Derivative estimates of semigroups and Riesz transforms on vector bundles». In: *Potential Anal.* 20.2 (2004), pp. 105–123.

[TW11] A. Thalmaier and F.-Y. Wang. «A stochastic approach to a priori estimates and Liouville theorems for harmonic maps». In: *Bull. Sci. Math.* 135.6-7 (2011), pp. 816–843.

[TW98] A. Thalmaier and F.-Y. Wang. «Gradient estimates for harmonic functions on regular domains in Riemannian manifolds». In: *J. Funct. Anal.* 155.1 (1998), pp. 109–124.

[Voi86] J. Voigt. «Absorption semigroups, their generators, and Schrödinger semigroups». In: *J. Funct. Anal.* 67.2 (1986), pp. 167–205.

[Wei80] J. Weidmann. *Linear operators in Hilbert spaces.* Vol. 68. Graduate Texts in Mathematics. Springer-Verlag, New York-Berlin, 1980, pp. xiii+402.